\documentclass[10pt]{article}
\usepackage{latexsym}
\usepackage[all]{xy}
\usepackage{amsmath}
\usepackage{amssymb}
\usepackage{amsfonts}

\newtheorem{thm}{Th\'eor\`eme}[section]
\newtheorem{prop}[thm]{Proposition}
\newtheorem{lem}[thm]{Lemme}

\newtheorem{df}[thm]{D\'efinition}
\newtheorem{cor}[thm]{Corollaire}

\begin{document}

\title{\textbf{Anneaux de Grothendieck des $n$-champs d'Artin}}
\bigskip
\bigskip

\author{\bigskip\\
Bertrand To\"en \\
\small{Institut de Math\'ematiques et de Mod\'elisation de Montpellier, UMR CNRS 5149}\\
\small{Universit\'e de Montpellier 2}\\
\small{Case Courrier 051}\\
\small{Place Eug\`ene Bataillon}\\
\small{34095 Montpellier Cedex}\\
\small{France}\\
\small{e-mail: btoen@math.univ-montp2.fr}}

\date{Novembre 2009}

\maketitle

\begin{abstract}
Nous introduisons un anneau de Grothendieck des champs d'Artin sup\'erieurs qui 
g\'en\'eralise la notion d'anneau de Grothendieck des vari\'et\'es. 
Nous montrons que cet anneau est non trivial en remarquant qu'il factorise l'invariant
\emph{nombre de points rationels sur un corps fini}. Dans une seconde partie nous
introduisons la notion de \emph{champs sp\'eciaux}, qui ont des groupes
d'homotopie $\pi_{i}$ affines, et unipotents pour $i>1$.
Notre th\'eor\`eme principal
affirme que le morphisme naturel de l'anneau de Grothendieck des vari\'et\'es
vers celui des champs sp\'eciaux induit un isomorphisme apr\`es 
inversion des classes de $\mathbb{A}^{1}$ et de $\mathbb{A}^{i}-\{0\}$ pour $i>0$. 
Nous d\'eduisons de ceci que de nombreux invariants num\'eriques (nombres de Hodge, 
caract\'eristique d'Euler motivique ou l-adique) s'\'etendent de fa\c{c}on unique
en des invariants de champs sp\'eciaux. En particulier, nous obtenons une version
de la formule des traces pour les champs d'Artin sp\'eciaux de type fini sur
un corps fini, qui identifie le nombre de points rationels \`a la trace du Frobenius
sur la carat\'eristique d'Euler l-adique \`a support compact. 
\end{abstract}

\medskip

\tableofcontents

\bigskip

\section{Introduction}

Les $n$-champs d'Artin sont des g\'en\'eralisations des 1-champs\footnote{Par la suite
le mot \emph{champ} signifiera un $n$-champ pour un certain $n$, et ainsi
nous utiliserons l'expression \emph{$1$-champ} pour faire r\'ef\'erence \`a la notion
usuelle telle qu'expos\'ee dans \cite{lm}.} alg\'ebriques au sens
de \cite{lm}, et qui apparaissent naturellement comme objets repr\'esentants certains
probl\`emes de modules qui ne sont (en g\'en\'eral) pas repr\'esentables par 
des 1-champs alg\'ebriques (voir \cite{to} pour une introduction plus d\'etaill\'ee).
Les bases d'une th\'eorie des $n$-champs d'Artin sont aujourd'hui \'etablies (voir \cite{si} et 
\cite[\S 2.1]{hagII}), et 
plusieurs exemples de tels objets ont \'et\'e construits (voir par exemple \cite{tv}). Il est alors naturel de 
se poser la question suivante : \emph{Que gagne-t-on \`a savoir qu'un probl\`eme
de modules est repr\'esentable par un $n$-champ d'Artin ?} Il existe de nombreuses approches
possibles \`a cette question tr\`es g\'en\'erale: on peut par exemple
se demander quels sont les invariants alg\'ebrico-g\'eom\'etriques d\'efinis pour les 
vari\'et\'es alg\'ebriques que l'on peut
\'etendre en des invariants de $n$-champs d'Artin. C'est ce que nous ferons dans ce travail, 
mais en se retreingant aux invariants additifs (e.g. caract\'eristiques d'Euler, nombre de points
sur un corps fini, polyn\^omes de Hodge \dots). 

Dans le cadre des vari\'et\'es alg\'ebriques (ou plus g\'en\'eralement des sch\'emas de type fini sur un anneau noeth\'erien $k$), 
il existe un invariant additif universel qui prend ses valeurs
dans un certain anneau de Grothendieck $K(\mathcal{V}(k))$, d\'efini comme le groupe ab\'elien
engendr\'e par les classes d'isomorphismes de vari\'et\'es et en imposant la
relation $[X]=[X-Y]+[Y]$ pour toute sous-vari\'et\'e ferm\'ee $Y \subset X$ (voir \cite{dl}, voir aussi 
Def. \ref{d7} o\`u nous utilisons une relation plus forte lorsque $k$ n'est pas de caract\'eristique nulle). 
Dans ce travail nous introduisons un anneau de Grothendieck $K(\mathcal{CH}^{ft}(k))$
des champs d'Artin sup\'erieurs de type fini sur $k$ (voir Def. \ref{d4}). La d\'efinition $K(\mathcal{CH}^{ft}(k))$ suit celle pour
les vari\'et\'es mais avec une relation suppl\'ementaire qui trivialise certaines 
fibrations localement triviales pour la topologie de Zariski (voir Def. \ref{d4} $(3)$). Il est \`a noter que cette nouvelle relation
est toujours satisfaite dans le cas des vari\'et\'es, mais est indispensable pour traiter
le cas des champs. Une premi\`ere observation est que l'anneau
$K(\mathcal{CH}^{ft}(k))$ est non-nul. Nous montrons cela lorsque 
$k=\mathbb{F}_{q}$ est un corps fini,  et en d\'efinissant un morphisme
d'anneaux 
$$\mu : K(\mathcal{CH}^{ft}(k)) \longrightarrow \mathbb{Q}$$
qui compte le nombre de points rationels en un sens convenable (i.e. en tenant compte de
l'existence des groupes d'homotopie sup\'erieurs des champs, voir Prop. \ref{p4}). 

Dans une seconde partie nous nous restreignons au cas 
des champs d'Artin \emph{sp\'eciaux}, qui par d\'efinition ont 
des groupes d'homotopie $\pi_{i}$ lin\'eaires et de plus unipotents pour $i>1$ (voir 
Def. \ref{d6}). Il se trouve que la classe des champs sp\'eciaux contient d\'ej\`a de nombreux
exemples int\'eressants, comme par exemple tous les $1$-champs d'Artin dont la diagonale est affine, 
ou encore les champs classifiants des structures lin\'eaires comme 
le champ $\mathcal{M}_{T}$ construit dans \cite{tv}. Les champs sp\'eciaux jouent aussi un r\^ole important
en th\'eorie de Hogde non-ab\'elienne o\`u ils sont appel\'es \emph{champs
tr\`es pr\'esentables} (voir \cite{si2,si3}). On consid\'erera alors l'anneau de Grothendieck
$K(\mathcal{CH}^{sp}(k))$ des champs sp\'eciaux d\'efini de mani\`ere analogue \`a l'anneau
$K(\mathcal{CH}^{tf}(k))$. Par construction, il existe un morphisme naturel
$$K(\mathcal{V}(k)) \longrightarrow K(\mathcal{CH}^{sp}(k))$$
induit par l'incusion des vari\'et\'es dans les champs d'Artin. Notre th\'eor\`eme principal est le suivant.

\begin{thm}\label{ti}
Soit $\mathbb{L}=[\mathbb{A}^{1}] \in K(\mathcal{V}(k))$ la classe de la droite affine,
et $\mathbf{1}:=[Spec\, k]$.  Notons
$K(\mathcal{V}(k))[\mathbb{L}^{-1},\{(\mathbb{L}^{i}-\mathbf{1})^{-1}\}_{i>0}]$ l'anneau localis\'e obtenu
en inversant $\mathbb{L}$ ainsi que tous les $\mathbb{L}^{i}-\mathbf{1}$ pour $i>0$. 
Alors, le morphisme induit
$$K(\mathcal{V}(k))[\mathbb{L}^{-1},\{(\mathbb{L}^{i}-\mathbf{1})^{-1}\}_{i>0}] \longrightarrow 
K(\mathcal{CH}^{sp}(k))$$
est un isomorphisme. Ainsi, tout invariant additif de vari\'et\'es \`a valeurs dans un anneau $A$
$$\phi : K(\mathcal{V}(k)) \longrightarrow A,$$
tel que $\phi(\mathbb{L})$ et les $\phi(\mathbb{L}^{i}-\mathbf{1})$ soient inversibles pour tout $i>0$,
s'\'etend de fa\c{c}on unique en un invariant additif de champs d'Artin sp\'eciaux
$$K(\mathcal{CH}^{sp}(k)) \longrightarrow A.$$
\end{thm}

Ce th\'eor\`eme r\'epond ainsi en partie 
\`a la question originale, de savoir quels sont les invariants de vari\'et\'es qui s'\'etendent
en des invariants de champs d'Artin sup\'erieurs. Sa d\'emonstration repose sur 
un d\'evissage des $n$-champs sp\'eciaux en des pieces \'el\'ementaires qui sont 
ou bien des sch\'emas ou bien  des champs
de la forme $K(\mathbb{G}_{a},n)$ avec $n>1$, ou bien des
1-champs quotients $[X/Gl_{r}]$ avec  $X$ un espace alg\'ebrique
muni d'une action du groupe $Gl_{r}$.  
Les pr\'erequis n\'ecessaires pour pouvoir effectuer ces d\'evissages
sont donn\'es dans les premiers paragraphes. 

En corollaire de ce th\'eor\`eme nous obtenons l'existence de nombreux 
invariants additifs de champs d'Artin sp\'eciaux, tel les nombres de Hodge ou encore
les caract\'eristiques d'Euler l-adiques et motiviques (voir \S 3.5). On d\'eduit aussi de ce th\'eor\`eme une
version de la formule des traces pour les champs sp\'eciaux, qui identifie le nombre
de points rationels sur un corps fini \`a la trace du Frobenius sur la caract\'eristique
d'Euler l-adique (voir Prop. \ref{p5}). \\

Deux remarques pour terminer cette introduction. Mise \`a part
l'invariant $\mu$ d\'efini dans la proposition \ref{d4} nous ne construisons pas dans ce travail
d'autres invariants d\'efinis sur $K(\mathcal{CH}^{tf}(k))$. Cependant, on peut
par exemple d\'efinir un invariant du type caract\'eristique d'Euler l-adique 
qui prolonge celui d\'efini pour les champs sp\'eciaux. Cette construction demande
l'existence d'un formalisme l-adique pour les champs sup\'erieurs qui 
n'est pour le moment pas disponible dans la litt\'erature (voir cependant \cite{be,lo} pour
le cas des $1$-champs, o\`u un tel formalisme est d\'evelopp\'e). 
Un travail non publi\'e (\cite{tvv}) devrait fournir un formalisme 
l-adique suffisemment sophistiqu\'e pour permettre la construction d'un tel 
invariant, et je renvoie donc \`a ce travail futur o\`u cette question sera probablement trait\'ee. 

Par ailleur, la relation $(3)$ de la d\'efinition \ref{d4} de l'anneau de Grothendieck
$K(\mathcal{CH}^{tf}(k))$ n'est vraiment justifi\'ee parceque l'on ne s'int\'eresse par la suite qu'aux champs d'Artin sp\'eciaux. 
Dans le cas g\'en\'eral il serait plus naturel d'utiliser la relation analogue, mais plus forte,
valable pour toute $F_{0}$-fibration Zariski localement triviale et pour tout champ
$F_{0}$. Il me semble que 
lorsque l'on se restreint aux champs sp\'eciaux cela donne un anneau de Grothendieck
isomorphe \`a celui d\'efini dans la d\'efinition \ref{d6}. Dans le cas g\'en\'eral cela donne un
anneau de Grothendieck probablement de taille plus raisonable que 
celui de notre d\'efinition \ref{d4}. Quoiqu'il en soit, l'analogue de notre th\'eor\`eme \ref{ti}
n'est de toute fa\c{c}on pas valable  pour des champs d'Artin g\'en\'eraux, et c'est la raison
pour laquelle je me suis content\'e des d\'efinitions \ref{d4} et \ref{d6}, certes moins naturelles mais 
certainement plus faciles \`a manipuler. \\

\bigskip

\textbf{Remerciements:} Je remercie M. Vaqui\'e et G. Vezzosi pour plusieurs discussions 
sur la notion d'invariants additifs de champs sup\'erieurs qui m'ont motiv\'e pour
\'ecrire ce texte. Je remercie aussi J. Sch\"urmann pour m'avoir signal\'e
une erreur dans une premi\`ere version de ce travail. 

\bigskip

\textbf{Conventions et notations:} Tout au long de ce travail $k$ sera un anneau
(associatif, commutatif et unitaire) noeth\'erien.

\bigskip

\section{Pr\'eliminaires sur les $n$-champs d'Artin}

Dans cette premi\`ere section, nous rappelons quelques d\'efinitions
sur les champs $n$-g\'eom\'etriques au sens de \cite[\S 2.1]{hagII}, et nous
en donnons quelques exemples \'el\'ementaires. Nous introduisons aussi 
la notion de gerbes et de gerbes totales, qui sont des analogues sup\'erieurs 
de la notion usuelle de gerbes. Un r\'esultat cl\'e 
pour la suite est que tout $n$-champ d'Artin fortement de pr\'esentation finie poss\`ede une stratification par
des gerbes totales. Nous utiliserons la d\'efinition de champs d'Artin bas\'ee sur la topologie 
fppf, et dont les atlas sont des atlas plats et de pr\'esentation finie. Cette notion, qui diff\`ere \`a priori  de la notion de champs d'Artin bas\'ee sur la topologie \'etale et atlas lisse
(telle que pr\'esent\'ee dans \cite{si,hagII}), lui est en r\'ealit\'e 
\'equivalente (voir \cite{fpqcet}). L'utilisation de la topologie plate, qui suppose implicitement les r\'esultats de \cite{fpqcet}, simplifiera les \'enonc\'es de d\'evissages en gerbes que nous allons  pr\'esenter de ce paragraphe.

\subsection{Rappel des d\'efinitions}

Rappelons que $k-Aff^{\sim,fppf}$ d\'esigne la cat\'egorie de
mod\`eles des pr\'efaisceaux simpliciaux sur le site des $k$-sch\'emas
affines munie de la topologie fid\`element plate et de pr\'esentation finie (\emph{fppf} pour faire court). 
On rappelle que les \'equivalences sont les
\'equivalences locales (e.g. les morphismes induisant des isomorphismes
sur tous les faisceaux d'homotopie, pour la topologie $fppf$), et que l'on utilise la structure
projective d\'ecrite dans \cite{bl} (voir aussi \cite{hagI}).
Nous passerons sous silence les 
histoires d'univers, et l'on renvoit \`a \cite{hagII} pour plus de d\'etails. La cat\'egorie
homotopique $St(k)$ est simplement appel\'ee \emph{la cat\'egorie
des champs}, et ses objets seront appel\'es
des \emph{champs}. Les ensembles de morphismes dans $St(k)$ seront not\'es $[-,-]$. De m\^eme un morphisme de champs sera toujours
un morphisme dans $St(k)$. Plus 
pr\'ecis\'ement, nous dirons qu'un champ $F$ 
est un \emph{$n$-champ}, si pour tout $i>n$, tout $X\in k-Aff$ et tout
$x\in F(X)$, le faisceau $fppf$ en groupes $\pi_{i}(F,x)$, associ\'e au pr\'efaisceau
$$\begin{array}{cccc}
\pi_{i}^{pr}(F,x) : & k-Aff/X & \longrightarrow & Gp \\
 & (f : Y \rightarrow X) & \mapsto & \pi_{i}(F(Y),f^{*}(x)),
 \end{array}$$
est trivial. Comme il est expliqu\'e dans \cite{hagI}, la sous-cat\'egorie
pleine de $St(k)$ form\'ee des $0$-champs est naturellement \'equivalente
\`a la cat\'egorie $Sh(k-Aff)$ des faisceaux en ensembles sur $k-Aff$. Le plongement de 
Yoneda $k-Aff \longrightarrow Sh(k-Aff)$ permet donc d'identifier la cat\'egorie
$k-Aff$ avec une sous cat\'egorie pleine de  $St(k)$. Par la suite, 
nous fairons cette identification de fa\c{c}on implicite.
 
Pour un objet $F \in St(k)$ nous noterons (en faisant un abus
de notation) $F(A)\in SEns$ la valeur
d'un remplacement fibrant de $F$ sur $A$. En d'autres termes on peut \'ecrire
$$F(A)\simeq Map(Spec\, A,F),$$
o\`u le $Map$ est calcul\'e dans la cat\'egorie de mod\`eles $k-Aff^{\sim,fppf}$
(noter que ce que nous notons $F(A)$ est not\'e $\mathbb{R}F(A)$
dans \cite{hagII}). 
Ceci signifie simplement que nous prenons soins de toujours effectuer
un remplacement fibrant (i.e. de passer au champ associ\'e) avant
d'\'evaluer en un point $A$. 
 
Nous dirons de mani\`ere \'equivalente que 
deux champs $F$ et $F'$ sont \emph{\'equivalents} ou bien qu'ils sont 
\emph{isomorphes dans} $St(k)$. \\

Avant de passer aux champs g\'eom\'etriques signalons que
la cat\'egorie $k-Aff^{\sim,fppf}$ est un \emph{topos de mod\`eles}, au sens
de \cite{hagI}. Ainsi, bien que la cat\'egorie des champs $St(k)$
ne soit pas un topos (par exemple elle ne poss\`ede pas
de limites), on dispose tout de m\^eme d'analogues homotopiques
du comportement bien connu des faisceaux. Par exemple, l'existence
de limites et de colimites homotopiques, les notions de monomorphismes
(ceux qui induisent des monomorphismes sur les faisceaux
$\pi_{0}$ et des isomorphismes sur tous les $\pi_{i}$ pour $i>0$) 
et d'\'epimorphismes (ceux qui induisent des \'epimorphismes de faisceaux
sur les faisceaux $\pi_{0}$), d'effectivit\'e de certains quotients \dots.
On renvoit \`a \cite{hagI} o\`u le lecteur pourra trouver
des pr\'ecision sur cette notion de topos de mod\`eles. 

Les produits fibr\'es homotopiques dans $St(k)$ seront not\'es
$F\times_{G}^{h}H$. \\

La d\'efinition des champs g\'eom\'etriques proc\`ede par la 
r\'ecurrence suivante. 
 
\begin{itemize}
 
 \item Un champ est $(-1)$-g\'eom\'etrique si c'est un sch\'ema affine. 
 
 \item Un morphisme de champs $F \longrightarrow G$ est 
 $(-1)$-repr\'esentable si pour tout $X\in k-Aff$, et tout
 morphisme de champs $X \longrightarrow G$, le 
 champ $F\times^{h}_{G}X$ est $(-1)$-g\'eom\'etrique. Un tel morphisme est 
 plat de pr\'esentation finie si de plus le morphisme induit
 $F\times^{g}_{G}X \longrightarrow X$ est un morphisme plat et de pr\'esentation finie entre sch\'emas affines.
 
\item Soit $n\geq 0$ un entier.
Supposons que pour tout $m<n$ la notion de champs $m$-g\'eom\'etriques soit d\'efinie, 
ainsi que  la notion de morphismes $m$-repr\'esentables et de morphismes $m$-repr\'esentables plats et de pr\'esentation finie.
 
  \begin{itemize}
  
 \item Un champ $F$ est $n$-g\'eom\'etrique s'il v\'erifie les deux conditions
 suivantes.
 
 \begin{itemize}

\item Le morphisme diagonal $F \longrightarrow F \times F$ est 
$(n-1)$-repr\'esentable.

\item Il existe des sch\'emas affines $X_{i}$ et un \'epimorphisme
de champs (i.e. dont le morphisme induit sur les faisceaux
$\pi_{0}$ soit un \'epimorphisme de faisceaux d'ensembles)
$$p : \coprod X_{i} \longrightarrow F,$$
tel que chaque morphisme $X_{i} \longrightarrow F$ 
(qui est $(n-1)$-repr\'esentable par la premi\`ere
condition sur la diagonale) soit plat et de pr\'esentation finie. Un tel 
morphisme $p$ est appel\'e un \emph{$n$-atlas}, ou simplement un 
atlas si l'on ne souhaite pas sp\'ecifier $n$.

 \end{itemize}
 
 \item Un morphisme $F \longrightarrow G$ est $n$-repr\'esentable si
 pour tout $X\in k-Aff$, et tout morphisme de champs $X \longrightarrow F$, 
 le champ $F\times_{G}^{h}X$ est $n$-g\'eom\'etrique. 
 
 \item Un morphisme $n$-repr\'esentable $F \longrightarrow G$ est lisse
 si pour tout $X\in k-Aff$ et tout morphisme de champs
 $X \longrightarrow F$, il existe un $n$-atlas
 $$\coprod X_{i} \longrightarrow F\times_{G}^{h}X,$$
tel que chacun des morphismes de sch\'emas 
$$X_{i} \longrightarrow X$$
soit plat et de pr\'esentation finie.
  
 \end{itemize}

\end{itemize}

Il est d\'emontr\'e dans \cite{hagII} que la notion de champs $n$-g\'eom\'etriques
se comporte comme on le souhaite. Par exemple,  la sous-cat\'egorie
des champs $n$-g\'eom\'etriques est stable par produits fibr\'es homotopiques, et 
\^etre $n$-g\'eom\'etrique est une propri\'et\'e locale pour la topologie
fppf. \\

Il est important de noter qu'un champ $F$ qui est $n$-g\'eom\'etrique
est automatiquement un $n$-champ (voir \cite{hagII}). En contre partie, un 
$n$-champ $F$ peut tout \`a fait \^etre $m$-g\'eom\'etrique pour
$m>n$, et ne pas \^etre $n$-g\'eom\'etrique. Ainsi, la complexit\'e
g\'eom\'etrique est toujours sup\'erieure \`a la complexit\'e homotopique, et
peut-\^etre strictement plus grande. Par exemple, 
un faisceau repr\'esentable par un sch\'ema sans aucune hypoth\`ese
de s\'eparation est un champ $1$-g\'eom\'etrique mais
pas $0$-g\'eom\'etrique en g\'en\'eral. Nous adopterons ainsi
la d\'efinition suivante.

\begin{df}\label{d1}
Un \emph{$n$-champ d'Artin} est un $n$-champ 
qui est $m$-g\'eom\'etrique pour un certain entier $m$. 
Un \emph{champ d'Artin} est un champ qui est un 
$n$-champ d'Artin pour un certain entier $n$.
\end{df} 

Le th\'eor\`eme principal de \cite{fpqcet} affirme que la notion pr\'ec\'edente de $n$-champs d'Artin coincide avec celle pr\'esent\'ee dans \cite{si,hagII}. On voit ainsi ais\'ement qu'un faisceau sur $k-Aff$ est repr\'esentable par un espace
alg\'ebrique si et seulement si c'est un 
$0$-champ d'Artin au sens de la d\'efinition ci-dessus. \\
 
Tout comme il est expliqu\'e dans \cite{hagII}, toute
propri\'et\'e \textbf{Q} de morphismes dans $k-Aff$ qui est locale
pour la topologie fppf, s'\'etent de fa\c{c}on naturelle en une propri\'et\'e de morphismes
entre $n$-champs d'Artin. Ainsi peut-on parler de morphismes
plats. Pour parler de morphismes lisses il est pr\'ef\'erable d'utiliser l'existence d'atlas lisses. Nous dirons ainsi qu'un morphisme de $n$-champs d'Artin $f : F \longrightarrow G$ est lisse, s'il existe 
des atlas lisses $U=\coprod_{i}U_{i} \rightarrow F$ et $V=\coprod_{i}V_{i} \rightarrow G$, et 
un diagramme commutatif
$$\xymatrix{
U \ar[r] \ar[d]_{g} & F \ar[d]^-{f} \\
V \ar[r] & G,}$$
avec $g$ lisse. \\

Rappelons aussi qu'un morphisme de champs $f : F \longrightarrow G$ est 
un monomorphisme, si le morphisme naturel 
$$F \longrightarrow F\times^{h}_{G}F$$
est un isomorphisme dans $St(k)$. De mani\`ere \'equivalente
$f$ induit un monomorphisme sur les faisceaux $\pi_{0}$ et 
des isomorphismes sur tous les faisceaux $\pi_{i}$ pour $i>0$ (et tout
choix de point de base). On d\'efinit alors les immersions ouvertes comme
\'etant les monomorphismes lisses. Les immersions ferm\'ees sont 
d\'efinies directement comme les morphismes $f : F \longrightarrow G$ 
tels que pour tout $X\in k-Aff$, et tout morphisme $X\longrightarrow G$, 
le  
morphisme $F\times_{G}^{h}X \longrightarrow X$ soit une immersion ferm\'ee de sch\'emas
affines. \\

Rappelons aussi les conditions de finitudes suivantes. 

\begin{itemize}
\item Un champ $F$ est 
\emph{quasi-compact} s'il existe un 
sch\'ema affine $X$ et un \'epimorphisme
de champs $X \longrightarrow F$. 

\item Un champ 
$0$-g\'eom\'etrique (i.e. un sch\'ema affine) est
\emph{fortement quasi-compact}.

\item Un morphisme $0$-repr\'esentable
$F \longrightarrow G$ est \emph{fortement quasi-compact}.

\item Par r\'ecurrence sur $n$, si $F$ est un champ 
 $n$-g\'eom\'etrique, nous
dirons que $F$ est \emph{fortement
quasi-compact} si les deux conditions suivantes 
sont satifaites.

\begin{itemize}
\item Le morphisme diagonal $F \longrightarrow F\times F$
est fortement quasi-compact. 

\item Le champ $F$ est quasi-compact.

\end{itemize}

\item Par r\'ecurrence sur $n$, nous dirons qu'un 
morphisme $n$-repr\'esentable $F \longrightarrow G$ 
est \emph{fortement quasi-compact} si pour tout sch\'ema affine
$X$ et tout morphisme $X \longrightarrow F$, le champ
$n$-g\'eom\'etrique $F\times_{G}^{h}X$ est 
fortement quasi-compact. 

\item Un  champ d'Artin $F$  est 
\emph{localement de pr\'esentation fini}  s'il existe un 
$n$-atlas $\coprod X_{i} \longrightarrow F$ o\`u chaque
$X_{i}$ est un sch\'ema affine 
de pr\'esentation finie sur $Spec\, k$. 

\item Un champ d'Artin $F$ est 
\emph{fortement de pr\'esentation fini} s'il est localement de pr\'esentation 
fini et fortement quasi-compact. 

\end{itemize}

Terminons cette premi\`ere section par
les touts premiers exemples de champ d'Artin, les
champs classifiants. Si $G$ est un sch\'ema en groupes 
ab\'eliens de pr\'esentation finie et plat sur
$Spec\, k$, alors on peut d\'efinir un pr\'efaisceau simplicial
$$\begin{array}{cccc}
K(G,n) : & k-Alg & \longrightarrow & SEns \\
 & A & \mapsto & K(G(A),n),
 \end{array}$$
o\`u comme d'habitude $K(H,n)$ est un ensemble simplicial
dont tous les groupes d'homotopie sont triviaux sauf 
$\pi_{n}(K(H,n))\simeq H$. Cet objet est un pr\'efaisceau simplicial
et sera donc consid\'er\'e comme objet dans 
$St(k)$. 
Il est facile de voir que 
$K(G,n)$ est un $n$-champ d'Artin, plat et fortement de
pr\'esentation finie sur $Spec\, k$. Le champ 
$K(G,n)$ est de plus lisse sur $Spec\, k$ d\`es que $n>0$, car 
le point global $* \longrightarrow K(G,n)$ est un atlas plat et de pr\'esentation: 
$K(G,n)$ est donc localement pour la topologie $fppf$ lisse sur $Spec\, k$, et est donc lisse sur 
$Spec\, k$.

Pour tout $k$-sch\'ema $X$, on dispose d'isomorphismes naturels
$$[X,K(G,n)] \simeq H_{fppf}^{n}(X,G).$$
Ainsi, d'apr\`es les r\'esultats de \cite{fpqcet}, on a 
$$[X,K(G,n)] \simeq H^{n}_{et}(X,G)$$
d\`es que $G$ est lisse. 

Enfin, lorsque $G$ n'est plus n\'ecessairement ab\'elien on dispose tout de m\^eme du $1$-champ 
classifiant $K(G,1)$, qui est lisse sur $Spec\, k$.

\subsection{Gerbes}

Soit $F\in St(k)$ un champ, et consid\'erons 
$\pi_{0}(F)$ le faisceau (pour la topologie fppf) 
associ\'e au pr\'efaisceau $X \mapsto \pi_{0}(F(X))$. 

\begin{df}\label{dgros}
Le \emph{faisceau des modules grossier} du champ $F$ est 
le faisceau $\pi_{0}(F)$. Nous le noterons $M(F)$. 
\end{df}

La notion de gerbe que l'on trouve dans la litt\'erature est 
vari\'ee, et elle diff\'ere suivant les contextes. Nous 
adopterons la terminologie suivante. 

\begin{df}\label{d2}
Un champ $F$ est une \emph{gerbe}
si les trois conditions suivantes sont satisfaites
\begin{enumerate}
\item Le champ $F$ est un champ d'Artin fortement de pr\'esentation finie. 
\item Le faisceau $M(F)$ est espace alg\'ebrique (i.e.
un $0$-champ d'Artin).
\item  Le morphisme naturel
$F \longrightarrow M(F)$ est plat. 
\end{enumerate}
\end{df}

Nous insistons sur le fait qu'une gerbe sera toujours pour nous un 
champ d'Artin fortement de pr\'esentation finie, ce qui n'est 
pas un terminologie tout \`a fait standard. \\

Nous rappelons que pour un champ
$F$, on d\'efinit le champ d'inertie $I_{F}$ (aussi appel\'e le champ 
des lacets) par
$$I_{F}:=\underline{Map}(S^{1},F)\simeq F\times_{F\times F}^{h}F,$$
o\`u $S^{1}$ est le champ associ\'e au pr\'efaisceau simplicial 
constant \'egal au cercle simplicial $S^{1}:=\Delta^{1}/\partial \Delta^{1}$, 
et o\`u $\underline{Map}$ d\'esigne le Hom interne de la cat\'egorie
homotopique des champs. 
Le choix d'un point de base $*\in S^{1}$ d\'efinit une projection naturelle
$$I_{F} \longrightarrow F,$$
qui peut aussi se voir comme l'une des deux projections
$$F\times_{F\times F}^{h}F \longrightarrow F.$$
Avec ces notations, le r\'esultat fondamental est le crit\`ere suivant, bien connu pour
les $1$-champs d'Artin (voir par exemple \cite{lm}).  

\begin{prop}\label{p1}
Un champ d'Artin $F$ est une gerbe si et seulement s'il est 
fortement de pr\'esentation finie et si de plus le morphisme naturel
$$I_{F} \longrightarrow F$$
est plat. 
\end{prop}

\textit{Preuve:} C'est la m\^eme que pour le cas $n=1$ trait\'e dans
\cite{lm}. R\'esumons-la bri\`evement. 

Commen\c{c}ons par voir que la condition est n\'ecessaire. Pour cela
on remarque que le fait que $I_{F} \longrightarrow F$ soit plat est 
une condition locale pour la topologie
plate de pr\'esentation finie sur $M(F)$. On peut donc supposer 
par changement de base que $M(F)=Spec\, k$, et m\^eme que la projection naturelle
$F \longrightarrow M(F)=Spec\, k$ poss\`ede une section, 
$x : Spec\, k \longrightarrow F$. Soit $U \longrightarrow F$ un morphisme fid\`element plat avec 
$U$ un sch\'ema affine (par exemple un atlas). Comme
$M(F)=Spec\, k$ on peut, quitte \`a prendre un recouvrement plat de $U$, supposer que
le morphisme $U \longrightarrow F$, se factorise par le point $x$
$$U \longrightarrow Spec\, k \longrightarrow F.$$
Alors, comme $U\longrightarrow F$ et $U\longrightarrow Spec\, k$ sont 
plats, on en d\'eduit que $x : Spec\, k \longrightarrow F$ est un morphisme plat. 
Ainsi, le champ des lacets en $x$
$$\Omega_{x}F:=Spec\, k\times_{F}^{h}Spec\, k$$
est un champ plat sur $Spec\, k$. Or, il existe un diagramme homotopiquement
cart\'esien de champs
$$\xymatrix{
\Omega_{x}F \ar[r] \ar[d] & Spec\, k \ar[d]^-{x} \\
I_{F} \ar[r] & F.}$$
Comme le morphisme $x : Spec\, k \longrightarrow F$
est fid\`element plat, ceci implique bien que le morphisme $I_{F} \longrightarrow F$
est plat.

Supponsons maintenant que $F$ soit un 
$n$-champ d'Artin fortement de pr\'esentation finie 
et que $I_{F} \longrightarrow F$ soit plat. 
Soit $X_{0} \longrightarrow F$ un atlas de $F$ 
(par quasi-compacit\'e on prendra
$X_{0}$ un sch\'ema affine), et 
consid\'erons
$$f : X_{1}:=X_{0}\times^{h}_{F}X_{0} \longrightarrow X_{0}\times X_{0}.$$
On note $R$ le sous-faisceau de $X_{0}\times X_{0}$
image du morphisme $f$. C'est le graphe d'une relation d'\'equivalence
sur $X_{0}$, et on a un isomorphisme de faisceaux 
$M(F)\simeq X_{0}/R$.  Pour montrer que $M(F)$ est un 
espace alg\'ebrique il nous faut donc montrer les deux assertions suivantes:

\begin{enumerate}
\item Le faisceau $R$ est repr\'esentable par un espace alg\'ebrique
localement de pr\'esentation finie.

\item Le morphisme $R \subset X_{0}\times X_{0} \longrightarrow X_{0}$
est plat. 
\end{enumerate}

Pour commencer, le morphisme diagonal $R \longrightarrow R\times R$
est un morphisme repr\'esentable affine car $R$ est un sous-faisceau
d'un sch\'ema affine. On consid\`ere de plus le morphisme
surjectif
$$X_{1} \longrightarrow R.$$
Pour tout sch\'ema affine $S$, et tout morphisme $S \longrightarrow R$, 
on dispose de deux morphismes induits $x : S \longrightarrow F$
et $y : S \longrightarrow F$, ou encore de deux points
$x$ et $y$ dans $F(S)$. Le produits fibr\'e homotopique
$X_{1}\times_{R}^{h}S$ s'identifie, en tant que champ sur $S$, 
au champ des chemins de $x$ \`a $y$
$$X_{1}\times_{R}^{h}S\simeq S\times_{F}^{h}S=:\Omega_{x,y}F
\longrightarrow S.$$
Le champ $\Omega_{x,y}F$ est localement non vide pour la topologie
fppf sur $S$ (car $x$ et $y$ proviennent d'un morphisme $S \longrightarrow R$). 
Ainsi, le champ $\Omega_{x,y}F$ est localement \'equivalent pour
la topologie fppf sur $S$ au champ $\Omega_{x}F$, des lacets
en $x$. Par hypoth\`ese sur $F$, le champ 
$\Omega_{x}F$ est plat de pr\'esentation finie sur $S$, et ainsi on voit que 
la projection 
$$X_{1}\times_{R}^{h}S \longrightarrow S$$
est un morphisme plat de pr\'esentation finie. Comme ceci est valable pour tout 
sch\'ema affine $S$ et tout morphisme $S \longrightarrow R$, on voit que
le morphisme de champs
$$X_{1} \longrightarrow R$$
est surjectif, plat et de pr\'esentation finie. En composant avec un atlas
pour $X_{1}$ on trouve donc un atlas pour $R$. 

Il nous reste \`a montrer que la premi\`ere projection 
$R \subset X_{0}\times X_{0} \longrightarrow X_{0}$ est un morphisme plat et de pr\'esentation finie. 
Mais, on dispose d'un diagramme commutatif de champs
$$\xymatrix{
X_{1} \ar[rd]^-{s} \ar[d]_-{p} & \\
R \ar[r]_-{q} & X_{0},}$$
et comme $s$ et $p$ sont plats de pr\'esentation finie et surjectifs il en est de m\^eme de 
$q$. Ceci termine la preuve que le faiscau 
$\pi_{0}$ est repr\'esentable par un espace alg\'ebrique de pr\'esentation finie, qui est le 
quotient $X_{0}/R$. De plus, on dispose d'un diagramme 
commutatif de champs
$$\xymatrix{
X_{0} \ar[rd]^-{a} \ar[d]_-{b} & \\
F \ar[r]_-{\pi} & M(F),}$$
avec $a$ et $b$ plats et surjectifs. Ainsi, le morphisme
$\pi$ est lui-m\^eme plat.  \hfill $\Box$ \\

\begin{cor}\label{c1}
Soit $F$ un champ d'Artin fortement de pr\'esentation fini. 
Alors $F$ admet une stratification par des sous-champs localement
ferm\'es qui sont des gerbes r\'eduites. 
En d'autres termes,
il existe une suite finie de sous-champs ferm\'es
$$\emptyset=F_{r+1} \subset F_{r} \subset \dots \subset 
F_{1}\subset F_{0}=F$$
tel que chaque champ $F_{i}-F_{i+1}$ soit une
gerbe r\'eduite.
\end{cor}

\textit{Preuve:} En posant $F_{1}:=F_{red}$ on se ram\`ene imm\'ediatement
au cas o\`u $F$ est r\'eduit.
Comme $k$ est noeth\'erien les champs
fortement de pr\'esentation finie sont noeth\'eriens (i.e. toute chaine
d\'ecroissante de sous-champs ferm\'es est stationaire). Ainsi, par r\'ecurrence
noeth\'erienne il nous suffit de voir que pour tout champ d'Artin fortement de pr\'esentation fini et r\'eduit
$F$, il existe un sous-champ ouvert $U\subset F$ tel que le champ
d'inertie 
$$I_{U}\simeq I_{F}\times_{F}^{h}U \longrightarrow U$$
soit plat. Mais ceci se d\'eduit ais\'ement de la platitude g\'en\'erique
pour les sch\'emas affines qui s'\'etend sans probl\`emes aux champs d'Artin.  \hfill $\Box$ \\

On peut \^etre plus pr\'ecis au sujet du corollaire \ref{c1}, en remarquant 
que tout stratification de $F$ par des sous-champs localement ferm\'es peut \^etre
rafin\'ee en une stratification par des gerbes r\'eduites. De plus, deux telles stratifications
par des gerbes r\'eduites peuvent \^etre rafin\'ees par une m\^eme troisi\`eme. 

\subsection{Repr\'esentabilit\'e des faisceaux d'homotopie}

Soit $F$ une gerbe r\'eduite. On rappelle que pour tout 
sch\'ema $X$ et tout morphisme $x : X \longrightarrow F$, on dispose de faisceaux
en groupes $\pi_{i}(F,x)$ sur $X$. 

\begin{prop}\label{p2}
Pour toute gerbe r\'eduite $F$, il existe un ouvert non-vide
$U$ de $F$, tel que pour tout sch\'ema $X$ affine, tout morphisme
$x : X \longrightarrow U$ et tout entier $i>0$, le faisceau
$\pi_{i}(F,x)$ soit repr\'esentable par espace alg\'ebrique en groupes
plat et de type fini sur $X$. 
\end{prop}

\textit{Preuve:} On raisone par r\'ecurrence sur un entier $n$ tel que
$F$ soit un $n$-champ d'Artin. 

Lorsque $n=1$, pour tout sch\'ema $X$ et tout morphisme
$x : X \longrightarrow F$ on a un isomorphisme de faisceaux sur $X$ 
$$\pi_{1}(F,x)\simeq I_{F}\times_{F}^{h}X.$$
Comme le morphisme $I_{F} \longrightarrow F$ est repr\'esentable de type fini, 
ceci implique que $\pi_{1}(F,x)$ est $n$-repr\'esentable (pour un certain $n$) par un espace
alg\'ebrique de type fini.
De plus, la proposition \ref{p1} implique que $\pi_{1}(F,x)$ est plat sur $X$. 
Ainsi, dans ce cas on peut prendre $U=F$. 

Supposons maintenant que $n>1$, et supposons la proposition d\'emontr\'ee pour
tous les $m$-champs d'Artin pour $m<n$. Notons $M=M(F)$ l'espace
de modules de $F$. Soit $X$ un sch\'ema affine, et 
$x : X \longrightarrow F$ un atlas lisse. Le champ $F$ \'etant un $n$-champ 
d'Artin, le faisceau $\pi_{n}(F,x)$ s'identifie au $n$-\`eme champ 
des lacets $\Omega_{x}^{(n)}F$ (on rappelle que l'on a par d\'efinition
$\Omega_{x}^{(n)}F=\Omega_{e}\Omega_{x}^{(n-1)}F$, o\`u $e$ est le point
distingu\'e de $\Omega_{x}^{(n-1)}F$ correspondant au lacet constant en $x$). 
Ceci implique que
$\pi_{n}(F,x)$ est repr\'esentable par un espace alg\'ebrique sur $X$. 
Soit $X'$ un ouvert non-vide de $X$ tel que $\pi_{n}(F,x)$ soit 
de plus plat sur $X$ (un tel ouvert existe car $F$, et donc $X$ est r\'eduit). 
L'image $F'$ du morphisme $X' \longrightarrow F$ est alors un sous-champ
ouvert qui est tel que pour tout sch\'ema $Y$ et tout morphisme
$y : Y \longrightarrow F'$, le faisceau $\pi_{n}(F,y)$ est
repr\'esentable par un espace alg\'ebrique plat et de type fini sur $Y$. 
On consid\`ere maintenant le morphisme de troncation
(voir \cite{hagI})
$$F \longrightarrow \tau_{\leq n-1}F.$$
Pour tout sch\'ema $X$ et tout morphisme 
$z : X\longrightarrow  \tau_{\leq n-1}F$, 
le morphisme $F\times_{\tau_{\leq n-1}F}^{h}X \longrightarrow X$
est localement (pour la topologie fppf) sur $X$ \'equivalent \`a un morphisme
de la forme $K(\pi_{n}(F,x),n) \longrightarrow X$, pour un 
point $x$ qui rel\`eve $z$. 
Or, le champ 
$K(\pi_{n}(F,x),n)$ est un $n$-champ d'Artin, et de plus le morphisme
$K(\pi_{n}(F,x),n) \longrightarrow X$ est un morphisme lisse. Ceci montre que 
le champ $F\times_{\tau_{\leq n-1}F}^{h}X$ est un $n$-champ d'Artin 
(voir \cite[Prop. 1.3.3.4]{hagII}), et donc que
 morphisme 
$$F \longrightarrow \tau_{\leq n-1}F$$
est $m$-repr\'esentable (pour un certain $m\geq n$) lisse, fortement de pr\'esentation finie, et surjectif. 
En composant avec un atlas pour $F$ on trouve donc un atlas pour $\tau_{\leq n-1}F$, ce qui montre que 
$\tau_{\leq n-1}F$ est un $(n-1)$-champ d'Artin.
Comme $F$ est de plus un champ d'Artin fortement de pr\'esentation finie, ceci implique que 
$\tau_{\leq n-1}F$
est aussi un champ d'Artin fortement de pr\'esentation finie. 

Le $(n-1)$-champ d'Artin $\tau_{\leq n-1}F$ est alors lui-m\^eme une gerbe. En effet, 
on a $M(\tau_{\leq n-1}F)=M(F)$.
De plus, 
comme les morphismes $F \longrightarrow \tau_{\leq n-1}F$ et 
$F \longrightarrow M(F)$ sont tous deux fid\`element plats, il en est de m\^eme du morphisme
$$\tau_{\leq n-1}F \longrightarrow M(\tau_{\leq n-1}F)=M(F).$$
Ainsi, par induction on sait qu'il existe un sous-champ ouvert $V$ 
de $\tau_{\leq n-1}F$ qui satisfait aux conditions de la proposition pour
$\tau_{\leq n-1}F$. L'image inverse de cet ouvert par le 
morphisme $F \longrightarrow \tau_{\leq n-1}F$ est le sous-champ ouvert
cherch\'e.
\hfill $\Box$ \\

Au vu de la proposition pr\'ec\'edente nous introduisons la d\'efinition
suivant.

\begin{df}\label{d2'}
Un champ $F$ est une \emph{gerbe totale} si c'est une gerbe, et si de plus
pour tout sch\'ema affine $X$, tout morphisme $x : X \longrightarrow F$ et tout
entier $i>0$, le faisceau $\pi_{i}(F,x)$ est repr\'esentable par un espace
alg\'ebrique plat sur $X$. 
\end{df}

De la proposition \ref{p2} et par r\'ecurrence noeth\'erienne on tire le corollaire
important suivant.

\begin{cor}\label{c1'}
Soit $F$ un champ d'Artin fortement de pr\'esentation fini. 
Alors $F$ admet une stratification par des sous-champs localement
ferm\'es qui sont des gerbes totales r\'eduites. 
En d'autres termes,
il existe une suite finie de sous-champs ferm\'es
$$\emptyset=F_{r+1} \subset F_{r} \subset \dots \subset 
F_{1}\subset F_{0}=F$$
tel que chaque champ $F_{i}-F_{i+1}$ soit une
gerbe totale r\'eduite.
\end{cor}

Le corollaire \ref{c1'} poss\`ede la cons\'equence important suivante.

\begin{cor}\label{c2}
Soit $F$ un champ d'Artin
fortement de pr\'esentation finie, $X$ un sch\'ema affine r\'eduit, et
et $x : X \longrightarrow F$ un morphisme. Alors, il existe un ouvert
non-vide $U$ de $X$, tel que
pour tout $i>0$, le faisceau
en groupes $\pi_{i}(F,x_{|U})$ (d\'efini sur $k-Aff/U$) soit repr\'esentable
par un espace alg\'ebrique en groupes plat et de type fini sur $U$.
En particulier, si $X=Spec\, K$ est le spectre d'un corps, alors
pour tout $i>0$ le faisceau $\pi_{i}(F,x)$ est repr\'esentable
par un sch\'ema en groupes de type fini sur $K$.
\end{cor}

\section{Anneaux de Grothendieck}

Nous d\'efinirrons dans cette section deux anneaux de Grothendieck, un premier pour les champs
d'Artin fortement de pr\'esentation finie et un second pour les champs d'Artin sp\'eciaux. 

\subsection{Anneau de Grothendieck des champs d'Artin fortement de pr\'esentation finie}

Nous commencerons par quelques d\'efinitions.

\begin{df}\label{d3}
\begin{enumerate}
\item 
Un morphisme $f : F \longrightarrow F'$ de champs d'Artin est 
une \emph{g-\'equivalence} si pour tout corps alg\'ebriquement clos $K$, le 
morphisme $F(K) \longrightarrow F'(K)$ est une \'equivalence.
\item Soit $F_{0}$ un champ.
Nous dirons qu'un morphisme de champs
$$F \longrightarrow F'$$
est une \emph{$F_{0}$-fibration Zariski localement triviale}
si pour tout sch\'ema affine $X$ et tout morphisme
$X \longrightarrow F'$, il existe un recouvrement 
Zariski $\{U_{i}\}$ de $X$, tel que chaque morphisme
$$F\times^{h}_{F'}U_{i} \longrightarrow U_{i}$$
soit isomorphe au-dessus de $U_{i}$ \`a la projection naturelle
$$F_{0}\times U_{i} \longrightarrow U_{i}.$$
\end{enumerate}
\end{df}

Notons maintenant $\mathcal{CH}^{tf}(k)$ la sous-cat\'egorie pleine de
$St(k)$ form\'ee des
champs d'Artin fortement de pr\'esentation finie
sur $k$. On consid\`ere le groupe ab\'elien libre $\mathbb{Z}[\mathcal{CH}^{tf}(k)]$ engendr\'e
par les classes d'isomorphismes de $\mathcal{CH}^{tf}(k)$. Comme il en est l'usage nous noterons 
$[F]\in \mathbb{Z}[\mathcal{CH}^{tf}(k)]$ la classe d'un champ d'Artin $F$. 

\begin{df}\label{d4}
Le \emph{groupe de Grothendieck des champs d'Artin} est le
quotient de $\mathbb{Z}[\mathcal{CH}^{tf}(k)]$ par les trois relations suivantes.
\begin{enumerate}
\item Pour $F$ et $F'$ dans $\mathcal{CH}^{tf}(k)$ on a
$$[F\coprod F']=[F]+[F'].$$
\item Pour toute $g$-\'equivalence $F \longrightarrow F'$ dans $\mathcal{CH}^{tf}(k)$, on a 
$$[F]=[F'].$$ 
\item Soit $F_{0}$ un champ qui est soit un sch\'ema affine, soit un 
champ de la forme $K(\mathbb{G}_{a},n)$ pour un entier $n>0$. 
Si $F \longrightarrow F'$ dans $\mathcal{CH}^{tf}(k)$ 
est une $F_{0}$-fibration Zariski localement triviale, alors on a
$$[F]=[F'\times F_{0}].$$
\end{enumerate}
Ce groupe sera not\'e $K(\mathcal{CH}^{tf}(k))$. 
\end{df}

On munit aussi $K(\mathcal{CH}^{tf}(k))$ d'une structure 
d'anneau commutatif en posant
$$[F].[F']:=[F\times F'],$$
et en v\'erifiant que ceci est compatible aux relations pr\'ec\'edentes. \\

Une remarque importante au sujet de la d\'efinition pr\'ec\'edente: toute
fibration en $K(\mathbb{G}_{a},n)$ qui est localement triviale pour
la topologie fppf est en r\'ealit\'e une $K(\mathbb{G}_{a},n)$-fibration
Zariski localement triviale. De fa\c{c}on plus pr\'ecise, si 
$f : F \longrightarrow F'$ est un morphisme de $n$-champs d'Artin
fortement de pr\'esentation finie, tel qu'il
existe un sch\'ema affine $X$ et 
un morphisme fi\`element plat de pr\'esentation finie $X \longrightarrow F'$
avec $F\times_{F'}X$ \'equivalent (comme champ sur $X$)
\`a $X\times K(\mathbb{G}_{a},n)$, alors $f$ est une
$K(\mathbb{G}_{a},n)$-fibration Zariski localement triviale. En effet, pour voir
cela il suffit de supposer que $F'$ est un sch\'ema affine $Y$.
Le champ
$A$ des auto-\'equivalences de $K(\mathbb{G}_{a},n)$ s'inscrit dans 
une diagramme homotopiquement cart\'esien
$$\xymatrix{
K(\mathbb{G}_{a},n) \ar[r] \ar[d] & A \ar[d] \\
Spec\, k \ar[r] & \mathbb{G}_{m}.} $$
Le champ $A$, comme champ en groupes, est m\^eme 
\'equivalent au produit semi-direct de $\mathbb{G}_{m}$
par $K(\mathbb{G}_{a},n)$. Ainsi, 
le champ classifiant des fibrations en $K(\mathbb{G}_{a},n)$ localement 
triviale pour la topologie fppf, 
qui n'est autre que $K(A,1)$, s'inscrit lui dans 
une diagramme homotopiquement cart\'sien
$$\xymatrix{
K(\mathbb{G}_{a},n+1) \ar[r] \ar[d] & K(A,1) \ar[d] \\
Spec\, k \ar[r] & K(\mathbb{G}_{m},1).} $$
Le morphisme $f : F \longrightarrow Y$ entre alors dans un diagramme homotopiquement 
cart\'esien
$$\xymatrix{
F \ar[r] \ar[d] & K(\mathbb{G}_{m},1) \ar[d] \\
Y \ar[r]_-{u} & K(A,1).}$$
Maintenant, en utilisant Hilbert 90, le morphisme
$u$ se factorise, localement pour la topologie de Zariski sur $Y$, 
en un morphisme $Y \longrightarrow K(\mathbb{G}_{a},n+1)$. 
Ce dernier morphisme correspond \`a $u\in H^{n+1}_{fppf}(Y,\mathbb{G}_{a})$, 
qui est donc nul localement pour la topologie Zariski sur $Y$ 
car $\mathbb{G}_{a}$ est le faisceau additif sous-jacent d'un 
faisceau coh\'erent. En conclusion on voit que le morphisme
$u : Y \longrightarrow K(A,1)$ est localement trivial pour la topologie
Zariski sur $Y$, ou de mani\`ere \'equivalente que 
$F\longrightarrow Y$ est une
$K(\mathbb{G}_{a},n)$-fibration Zariski localement triviale. Ce fait remarquable
sera utilis\'e implicitement par la suite.\\

Lorsque $k$ est de caract\'eristique nulle la notion de g-\'equivalence se simplifie, et le groupe
de Grothendieck $K(\mathcal{CH}^{tf}(k))$ poss\`ede alors une pr\'esentation 
plus proche de la d\'efinition usuelle du groupe de Grothendieck des vari\'et\'es (tel que d\'efinis dans \cite{dl}). 

\begin{prop}\label{p3}
Supposons que $k$ soit de caract\'eristique nulle. Alors le groupe
$K(\mathcal{CH}^{tf}(k))$ est isomorphe au quotient de
$\mathbb{Z}[\mathcal{CH}^{tf}(k)]$ par les trois relations suivantes.
\begin{enumerate}
\item Pour $F$ et $F'$ dans $\mathcal{CH}^{tf}(k)$ on a
$$[F\coprod F']=[F]+[F'].$$
\item Pour tout $F$ dans $\mathcal{CH}^{tf}(k)$, et 
tout sous-champ ferm\'e $F_{0} \subset F$, d'ouvert compl\'ementaire
$F-F_{0}$, on a
$$[F]=[F_{0}]+[F-F_{0}].$$ 
\item  Soit $F_{0}$ un champ qui est soit un sch\'ema affine, soit un 
champ de la forme $K(\mathbb{G}_{a},n)$ pour un entier $n>0$. Si $F \longrightarrow F'$ dans $\mathcal{CH}^{tf}(k)$ 
est une $F_{0}$-fibration Zariski localement triviale, alors on a
$$[F]=[F'\times F_{0}].$$
\end{enumerate} 
\end{prop}

\textit{Preuve:} Notons $K'(\mathcal{CH}^{tf}(k))$ le groupe d\'efini par les relations
de la proposition. On dispose d'un morphisme \'evident
$$\mathbb{Z}[\mathcal{CH}^{tf}(k)] \longrightarrow K(\mathcal{CH}^{tf}(k)).$$
De plus, pour $F$ dans $\mathcal{CH}^{tf}(k)$, et 
tout sous-champ ferm\'e $F_{0} \subset F$, le morphisme
$$F_{0} \coprod F-F_{0} \longrightarrow F$$
est une g-\'equivalence. Ceci implique que le morphisme ci-dessus induit un morphisme
$$K'(\mathcal{CH}^{tf}(k)) \longrightarrow K(\mathcal{CH}^{tf}(k)).$$
Pour montrer que cet un isomophisme, il suffit de montrer que le morphisme
\'evident
$$\mathbb{Z}[\mathcal{CH}^{tf}(k)] \longrightarrow K'(\mathcal{CH}^{tf}(k))$$
se factorise par $K(\mathcal{CH}^{tf}(k))$, car cela permet de construire un inverse
$$K(\mathcal{CH}^{tf}(k)) \longrightarrow K'(\mathcal{CH}^{tf}(k)).$$
Pour cela il suffit de voir que toute g-\'equivalence 
$f : F \longrightarrow F'$ dans $\mathcal{CH}^{tf}(k)$, on ait
$[F]=[F']$ dans $K'(\mathcal{CH}^{tf}(k))$. 

\begin{lem}\label{l1}
Soit $f : F \longrightarrow F'$ une g-\'equivalence entre champs d'Artin r\'eduits dans
$\mathcal{CH}^{tf}(k)$. Alors il existe un sous-champ ouvert non-vide 
$V \subset F'$, tel que le morphisme induit
$$F\times^{h}_{F'}V \longrightarrow V$$
soit un isomorphisme dans $\mathcal{CH}^{tf}(k)$.
\end{lem}

\textit{Preuve du lemme:} Soit $X' \longrightarrow F'$ un atlas avec $X'$ un sch\'ema affine r\'eduit. 
On consid\`ere le diagramme homotopiquement cart\'esien
$$\xymatrix{
X \ar[d] \ar[r] & X' \ar[d] \\
F \ar[r] & F'.}$$
Par hypoth\`ese sur $f$, pour tout corps alg\'ebriquement clos 
$K$, le morphisme $X(K) \longrightarrow X'(K)$ est un \'equivalence, et en 
particulier $X(K)$ est homotopiquement discret (i.e. \'equivalent \`a un ensemble). 
Montrons que cela implique que $X$ est un espace alg\'ebrique. 
Pour cela nous montrerons plus g\'en\'eralement qu'un champ
$n$-g\'eom\'etrique $Y$, fortement de pr\'esentation finie, et 
tel que $Y(K)$ soit (homotopiquement) discret pour tout corps alg\'ebriquement
clos $K$, est un espace alg\'ebrique. Ceci se d\'emontre par induction sur $n$
comme suit. Soit $Y_{*}$ un groupoide de Segal repr\'esentant 
$Y$, avec $Y_{0}$ un sch\'ema, et $Y_{1}$ un champ $(n-1)$-g\'eom\'etrique lisse
sur $Y_{0}$ (voir \cite{hagII}). Comme $K$ est alg\'ebriquement clos, le morphisme
$$Y_{0}(K) \longrightarrow Y(K)$$
est surjectif \`a homotopie pr\`es, et donc 
$Y(K)$ s'\'ecrit comme la r\'ealisation g\'eom\'etrique de $Y_{*}(K)$
$$Y(K)\simeq  |Y_{*}(K)|.$$
Comme l'ensemble simplicial $Y_{0}(K)$ est discret, 
les fibres homotopiques du morphisme
$$Y_{1}(K) \longrightarrow Y_{0}(K)\times Y_{0}(K)$$
sont \'equivalentes aux espaces de chemins de $Y(K)$, et sont 
donc ou bien vides ou bien contractiles. En d'autres termes,
$Y_{1}(K)$ est lui-m\^eme discret et  
le morphisme
$$Y_{1}(K) \longrightarrow Y_{0}(K)\times Y_{0}(K)$$
est un monomorphisme. Ceci \'etant valide pour tout corps alg\'ebriquement clos $K$, par r\'ecurrence on d\'eduit que 
$Y_{1}$ est un espace alg\'ebrique r\'eduit. Les fibres g\'eom\'etriques
du morphisme 
$$Y_{1} \longrightarrow Y_{0}\times Y_{0}$$
sont ou bien vide, ou bien des torseurs sous des sch\'emas en groupes. Comme
$k$ est de caract\'eristique nulle ces fibres g\'eom\'etriques sont donc r\'eduites, et 
m\^eme des spectres de corps car elles  poss\`edent au plus un point ferm\'e. 
Ainsi, on voit que le morphisme $Y_{1} \longrightarrow Y_{0}\times Y_{0}$ est non ramifi\'e. 
Comme il est de plus injectif sur les points g\'eom\'etriques il s'agit d'un monomorphisme. 
Ainsi, le champ $Y$, qui est le quotient du groupoide
de Segal $Y_{*}$, est un espace alg\'ebrique car quotient de $Y_{0}$
par la relation d'\'equivalence lisse $Y_{1} \longrightarrow Y_{0}\times Y_{0}$. 

Revenons \`a la preuve du lemme et \`a notre diagramme homotopiquement cart\'esien
 $$\xymatrix{
X \ar[d] \ar[r] & X' \ar[d] \\
F \ar[r] & F'.}$$
Nous savons maintenant que $X$ est un espace alg\'ebrique r\'eduit, et aussi que 
le morpisme $X \longrightarrow X'$ induit un bijection 
$X(K) \longrightarrow X'(K)$ pour tout corps alg\'ebriquement clos $K$. 
Comme nous sommes en caract\'eristique nulle, on sait qu'il existe un 
ouvert non-vide $U$ de $X'$, tel que le morphisme induit
$$X\times^{h}_{X}U \longrightarrow U$$
soit un isomorphisme. L'image $V$ de cet ouvert dans $F'$ v\'erifie clairement 
les conditions du lemme.  \hfill $\Box$ \\

Revenons \`a la preuve de la proposition. Il
il suffit de voir que toute g-\'equivalence 
$f : F \longrightarrow F'$ dans $\mathcal{CH}^{tf}(k)$, on ait
$[F]=[F']$ dans $K'(\mathcal{CH}^{tf}(k))$. Tout d'abord, 
comme $F_{red} : F_{red} \longrightarrow F'_{red}$ est 
encore une g-\'equivalence, et que 
$[F_{red}]=[F]$, $[F'_{red}]=[F']$ dans $K'(\mathcal{CH}^{tf}(k))$, on peut supposer
$F$ et $F'$ r\'eduits. On remarque aussi que la retriction de 
$f$ \`a un sous-champ localement ferm\'e de $F'$ reste une g-\'equivalence. 
Des applications successives du lemme \ref{l1}, impliquent alors que l'on puisse trouver
des suites d\'ecroissantes de ferm\'es r\'eduits 
$$F_{r+1}=\emptyset \subset F_{r} \subset \dots F_{1} \subset F_{0}=F$$
$$F'_{r+1}=\emptyset \subset F'_{r} \subset \dots F'_{1} \subset F'_{0}=F',$$
avec $F_{i}=(F'_{i}\times^{h}_{F'}F)_{red}$, et tel que les morphismes induits
$$F_{i}-F_{i-1} \longrightarrow F'_{i}-F'_{i-1}$$
soient des isomorphismes pour tout $i$. On a alors
$$[F]=\Sigma_{i}[F_{i}-F_{i-1}]=\Sigma_{i}[F'_{i}-F'_{i-1}]=[F']$$
dans $K'(\mathcal{CH}^{tf}(k))$. 
\hfill $\Box$ \\

\subsection{Nombre de points rationels}

Avant de passer aux cas des champs sp\'eciaux et \`a notre th\'eor\`eme
principal, signalons le fait suivant qui montre que le groupe
de Grothendieck $K'(\mathcal{CH}^{tf}(k))$ est non trivial et 
qu'il mesure bien la \emph{taille} des champs d'Artin. Pour cela nous
allons montrer qu'il 
factorise l'invariant \emph{nombre de points rationels} sur un 
corps fini. 

\begin{prop}\label{p4}
Supposons que $k=\mathbb{F}_{q}$ soit un corps fini. 
\begin{enumerate}
\item Pour tout champ d'Artin fortement de pr\'esentation finie
$F \in \mathcal{CH}^{tf}(k)$, l'ensemble $\pi_{0}(F(k))$ est fini. 
De plus, pour tout entier $i>0$ et tout $x \in \pi_{0}(F(k))$, le
groupe $\pi_{i}(F(k),x)$ est fini. 
\item L'application 
$$\mu : \mathbb{Z}[\mathcal{CH}^{tf}(k)] \longrightarrow \mathbb{Q}$$
d\'efinie par
$$\mu([F]):=\Sigma_{x\in \pi_{0}(F(k))}\prod_{i>0}|\pi_{i}(F(k),x)|^{(-1)^{i}},$$
(pour $F \in \mathcal{CH}^{tf}(k)$ et o\`u $|A|$ d\'esigne 
le cardinal d'un ensemble fini $A$) se factorise en un morphisme
d'anneaux
$$\mu : K(\mathcal{CH}^{tf}(k)) \longrightarrow \mathbb{Q}.$$
\end{enumerate}
\end{prop}

\textit{Preuve:} $(1)$ Il suffit de montrer que pour tout
champ d'Artin fortement de pr\'esentation finie $F$ l'ensemble
$\pi_{0}(F(k))$ est fini. En effet, pour un point 
$x\in F(k)$, on a 
$$\pi_{i}(F(k),x)\simeq \pi_{0}(\Omega_{x}^{(i)}F(k)),$$
o\`u $\Omega_{x}^{(i)}F$ est le $i$-\`eme champ de lacets 
au point $x$. Comme $\Omega_{x}^{(i)}F$ est encore un champ 
d'Artin fortement de pr\'esentation finie le r\'esultat pour
$\pi_{0}$ implique celui pour les $\pi_{i}$.

Soit donc $F \in \mathcal{CH}^{tf}(k)$ et montrons que
$\pi_{0}(F(k))$ est fini. 
En appliquant la proposition \ref{p1} on voit que l'on se ram\`ene au cas o\`u 
$F$ est une gerbe. Soit $F \longrightarrow M(F)=:M$ la projection de
$F$ sur son espace de modules. On dispose d'un morphisme induit
$$F(k) \longrightarrow M(k),$$
et donc pour tout point $x \in M(k)$ d'un suite exacte en homotopie
$$\xymatrix{\pi_{0}(F_{0}(k)) \ar[r] & \pi_{0}(F(k)) \ar[r] & M(k), }$$
o\`u $F_{0}$ est la fibre homotopique de $F \longrightarrow M$ en $x$.
Comme $M(k)$ est fini, 
ceci montre que $\pi_{0}(F(k))$ est fini si 
$\pi_{0}(F_{0}(k))$ l'est pour tout choix de $x \in M(k)$.
On se ram\`ene ainsi au cas o\`u 
$M=Spec\, k$, ou encore o\`u $M(F)=*$. 
On peut aussi clairement supposer que 
$F(k)$ est non vide, et on choisit donc
$x \in F(k)$.

En appliquant le corollaire \ref{c2} on voit alors que les 
faisceaux $\pi_{i}(F,x)$ sont repr\'esentables par
des espaces alg\'ebriques de type fini sur $k$.
Ceci a pour cons\'equence
que le champ $1$-tronqu\'e $\tau_{\leq 1}F$, qui est 
de la forme $K(\pi_{1}(F,x),1)$, est un champ d'Artin. 
Notons $H$ le sch\'ema en groupes $\pi_{1}(F,x)$, et 
$H_{0}$ sa composante neutre et r\'eduite.
Par le th\'eor\`eme de Lang, on sait que l'ensemble simplicial
$K(H_{0},1)(k)$ est connexe, et donc que le morphisme
$$\pi_{0}(K(H,1)(k)) \longrightarrow \pi_{0}(K(H/H_{0},1)(k))$$
est injectif. De plus, comme le sch\'ema en groupes 
$H/H_{0}=:K$ est fini, on peut trouver un plongement 
$K \hookrightarrow Gl_{n}$. On a alors
$$K(K,1)\simeq [X/Gl_{n}]$$
o\`u $X=Gl_{n}/K$. Ainsi, le morphisme $X \longrightarrow K(K,1)$
est un $Gl_{n}$-torseur, et donc le morphisme 
$$X(k) \longrightarrow \pi_{0}(K(K,1)(k))$$
est surjectif. En particulier le membre de droite est fini, et on conclut que 
l'ensemble
$\pi_{0}(\tau_{\leq 1}F(k))$ est lui-m\^eme fini. 

Pour conclure il nous reste a montrer que le morphisme naturel
$$F \longrightarrow \tau_{\leq 1}F$$
induit une injection
$$\pi_{0}(F(k)) \hookrightarrow \pi_{0}(\tau_{\leq 1}F(k)).$$
Pour cela il suffit de montrer que pour tout 
$x \in  \pi_{0}(F(k))$, l'ensemble simplicial 
$F_{>1}(k)$, fibre homotopique du morphisme pr\'ec\'edent en l'image $y$ de $x$, est 
connexe. Soit donc $F_{>1}$ le champ d\'efini par
$$F_{>1}:=F\times^{h}_{\tau_{\leq 1}F}*,$$
fibre en $y$ du morphisme $F \longrightarrow \tau_{\\eq 1}F$.
Le point $x$ d\'efinit un morphisme $z : * \longrightarrow F_{>1}$. 

\begin{lem}\label{l2}
Soit $F$ un $n$-champ d'Artin fortement de pr\'esentation finie muni d'un 
point global $x \in \pi_{0}(F(k))$. 
On suppose que $M(F)\simeq *$
ainsi que $\pi_{1}(F,x)\simeq 0$. Alors l'ensemble simplicial 
$F(k)$ est connexe. 
\end{lem}

\textit{Preuve du lemme:} Remarquons que le corollaire \ref{c1'} implique
que $F$ est un gerbe totale. En particulier tous les tronqu\'es
$\tau_{\leq m}F$ sont des champs d'Artin fortement de pr\'esentation finie.

Montrons cette proposition pour les $n$-champs
avec $n>1$ fix\'e (pour $n<2$ il n'y a rien \`a d\'emontrer).
On raisonne alors par r\'ecurrence sur l'entier $n$. 
Commen\c{c}ons par
le cas $n=2$. Par hypoth\`ese le champ 
$F$ est de la forme $K(H,2)$, pour $H$ un 
sch\'ema en groupes ab\'eliens de type fini. On a alors
$$\pi_{0}(F(k),x)\simeq H^{2}_{fppf}(Spec\, k,H)=0$$
car $k$ est de dimension cohomologique $1$ (voir \cite{dg} ou encore \cite[Cor. 3.2 (4)]{fpqcet}). 

Supposons maintenant que $n>2$, et aussi que 
le lemme soit d\'emontr\'e pour tous les $m$-champs $F$
avec $m<n$.  
On consid\`ere la troncation
$$F \longrightarrow \tau_{\leq n-1}F,$$
et sa fibre homotopique au point $x$. 
$$F_{n} \longrightarrow F.$$
Par hypoth\`ese, on a 
$F_{n}\simeq K(H,n)$
pour un certain sch\'ema en groupes ab\'eliens $H$ de type fini, et donc
comme ci-dessus
$$\pi_{0}(F_{n}(k),x)\simeq H^{n}_{fppf}(Spec\, k,H)=0.$$
De plus, par induction on a aussi $\pi_{0}(\tau_{\leq n-1}F(k))=*$. 
La suite longue d'homotopie associ\'ee \`a la fibration
$$\xymatrix{
F_{n}(k) \ar[r] & F(k) \ar[r] & \tau_{\leq n-1}F(k),}$$
et l'hypoth\`ese de r\'ecurrence
montrent alors que $F(k)$ est connexe. \hfill $\Box$ \\

Le lemme permet de conclure que le morphisme
$$\pi_{0}(F(k)) \hookrightarrow \pi_{0}(\tau_{\leq 1}F(k))$$
est injectif, et donc que 
$\pi_{0}(F(k))$ est fini par ce que l'on avait vu. Ceci termine le point $(1)$ de la proposition. \\

$(2)$ Il faut v\'erifier que $\mu$ est compatible aux trois
relations d\'efinissant $K(\mathcal{CH}^{tf})$. La premi\`ere est 
\'evidente.
Soit $F \longrightarrow F'$ est un g-\'equivalence. Elle
induit une \'equivalence 
$$F(\overline{k}) \longrightarrow F'(\overline{k})$$
qui est de plus $Gal(\overline{k}/k)$-\'equivariante (au sens profini). 
En passant au points fixes homotopiques (au sens profini) on trouve donc
que le morphisme induit
$$F(k)\simeq F(\overline{k})^{hGal(\overline{k}/k)} \longrightarrow 
F'(k)\simeq F'(\overline{k})^{hGal(\overline{k}/k)}$$
est un \'equivalence. Ceci implique clairement que 
$\mu([F])=\mu([F'])$, et donc que la seconde relation est v\'erifi\'ee. 

Il nous reste \`a montrer que $\mu$ est compatible avec la troisi\`eme relation. 
Soit $F \longrightarrow F'$  une $F_{0}$-fibration Zariski localement triviale. 
Il nous faut montrer que $\mu([F])=\mu([F']).\mu([F_{0}])$. 
Le morphisme
$$\pi_{0}(F(k)) \longrightarrow \pi_{0}(F'(k))$$
est clairement surjectif (\`a cause de la Zariski locale trivialit\'e), 
et ses fibres homotopiques sont  toutes
\'equivalentes \`a $F_{0}(k)$. La suite longue en homotopie
associ\'ee au morphisme $F(k) \longrightarrow F'(k)$
implique alors facilement que l'on a 
$$\mu([F])=\mu([F']).\mu([F_{0}])$$
(voir \cite{to2} pour plus de d\'etails).  
\hfill $\Box$ \\

\subsection{Anneau de Grothendieck des champs d'Artin sp\'eciaux}

Rappelons d'apr\`es le corollaire \ref{c2} que pour tout champ d'Artin
fortement de pr\'esentation finie $F$, tout corps $K$, tout
morphisme $x : Spec\, K \longrightarrow F$, le faisceau
$\pi_{i}(F,x)$ est repr\'esentable par un sch\'ema en groupes 
de type fini sur $K$.

\begin{df}\label{d5}
Un champ d'Artin $F$ est \emph{sp\'ecial}
s'il est fortement de pr\'esentation finie et s'il satisfait aux
deux conditions suivantes.
\begin{enumerate}
\item
Pour tout corps $K$, tout
point $x \in F(K)$, et tout $i>0$ 
le sch\'ema en groupes $\pi_{i}(F,x)$ est affine.
\item  Pour tout corps $K$, tout
point $x \in F(K)$, et tout $i>1$ 
le sch\'ema en groupes $\pi_{i}(F,x)$ est unipotent.
\end{enumerate} 
La sous-cat\'egorie pleine de 
$\mathcal{CH}^{tf}(k)$ form\'ee des champs d'Artin sp\'eciaux
sera not\'ee 
$\mathcal{CH}^{sp}(k)$.
\end{df}

On note $\mathbb{Z}[\mathcal{CH}^{sp}(k)]$ le groupe ab\'elien libre
engendr\'e par les classes d'isomorphismes de 
$\mathcal{CH}^{sp}(k)$. 
On d\'efinit alors le groupe de Grothendieck des champs
d'Artin sp\'eciaux de mani\`ere analogue \`a celle de la d\'efinition \ref{d4}. 

\begin{df}\label{d6}
Le \emph{groupe de Grothendieck des champs d'Artin sp\'eciaux} est le
quotient de $\mathbb{Z}[\mathcal{CH}^{sp}(k)]$ par les trois relations suivantes.
\begin{enumerate}
\item Pour $F$ et $F'$ dans $\mathcal{CH}^{sp}(k)$ on a
$$[F\coprod F']=[F]+[F'].$$
\item Pour toute $g$-\'equivalence $F \longrightarrow F'$ dans $\mathcal{CH}^{sp}(k)$, on a 
$$[F]=[F'].$$ 
\item Soit $F_{0}$ un champ qui est soit un sch\'ema affine, soit un 
champ de la forme $K(\mathbb{G}_{a},n)$ pour un entier $n>0$. 
Si $F \longrightarrow F'$ dans $\mathcal{CH}^{sp}(k)$ 
est une $F_{0}$-fibration Zariski localement triviale, alors on a
$$[F]=[F'\times F_{0}].$$
\end{enumerate}
Ce groupe sera not\'e $K(\mathcal{CH}^{sp}(k))$. 
\end{df}

Comme pour le cas des champs d'Artin fortement de pr\'esentation finie, 
le produit direct fait de $K(\mathcal{CH}^{sp}(k))$ un anneau commutatif. \\

Bien entendu on dispose d'un morphisme naturel
$$K(\mathcal{CH}^{sp}(k)) \longrightarrow 
K(\mathcal{CH}^{tf}(k)).$$
De m\^eme, soit $\mathcal{V}(k)$ la cat\'egorie des vari\'et\'es
sur $k$ (i.e. des sch\'emas de type fini sur $k$), et 
$\mathbb{Z}[\mathcal{V}(k)]$ le groupe ab\'elien libre engendr\'e par
ses classes d'isomorphismes. 
On d\'efinit alors  le groupe de Grothendieck des vari\'et\'es de la fa\c{c}on suivante.

\begin{df}\label{d7}
Le \emph{groupe de Grothendieck des vari\'et\'es} est le
quotient de $\mathbb{Z}[\mathcal{V}(k)]$ par les deux relations suivantes.
\begin{enumerate}
\item Pour $X$ et $Y$ dans $\mathcal{V}(k)$ on a
$$[X\coprod Y]=[X]+[Y].$$
\item Pour toute $g$-\'equivalence $X \longrightarrow Y$ dans $\mathcal{V}(k)$, on a 
$$[X]=[Y].$$ 
\end{enumerate}
Ce groupe sera not\'e $K(\mathcal{V}(k))$. 
\end{df}

Le produit direct fait de $K(\mathcal{V}(k))$ un anneau commutatif. \\

Le lecteur v\'erifiera \`a l'aide du lemme \ref{l1} que
lorsque $k$ est de caract\'eristique nulle $K(\mathcal{V}(k))$
est isomorphe au groupe de Grothendieck des vari\'et\'es que 
l'on rencontre habituellement dans la litt\'erature, d\'efini par l'unique
relation
$$[X]=[X-Y]+[Y]$$
pour tout sous-sch\'ema ferm\'e $Y \subset X$ (voir \cite{dl}).
Ceci dit, en caract\'eristique positive cette unique relation n'est 
pas suffisante pour retrouver le groupe
$K(\mathcal{V}(k))$ d\'efinit ci-dessus. 

Noter aussi que l'analogue de la relation $(3)$ des d\'efinitions 
\ref{d4} et \ref{d6} n'apparait pas dans la d\'efinition 
de $K(\mathcal{V}(k))$. Ceci s'explique par le fait que 
dans $K(\mathcal{V}(k))$, si $X \longrightarrow Y$ est une fibration 
Zariski triviale de fibre $Z$ alors on a toujours
$[X]=[Y].[Z]$. En effet, par stratification il suffit de montrer
cela g\'en\'eriquement sur $Y$, mais dans ce cas on peut supposer
que la fibration est triviale et c'est alors imm\'ediat. \\

L'inclusion des vari\'et\'es dans les champs d'Artin sp\'eciaux 
induit un morphisme d'anneaux bien d\'efini
$$j : K(\mathcal{V}(k)) \longrightarrow K(\mathcal{CH}^{sp}(k)).$$
Nous noterons $\mathbb{L}:=[\mathbb{A}^{1}]$, comme
\'el\'ement de $K(\mathcal{V}(k))$ ou encore de 
$K(\mathcal{CH}^{sp}(k))$ (noter que $j(\mathbb{L})=\mathbb{L}$).
De m\^eme, nous noterons $\mathbf{1}:=[Spec\, k]$, comme
\'el\'ement de $K(\mathcal{V}(k))$ ou encore de 
$K(\mathcal{CH}^{sp}(k))$. On remarque que 
le morphisme 
$$Spec\, k \longrightarrow K(Gl_{n},1)$$
est une $Gl_{n}$-fibration Zariksi localement triviale, et donc dans 
$K(\mathcal{CH}^{sp}(k))$ on a l'\'egalit\'e
$$[K(Gl_{n},1)].[Gl_{n}]=[Spec\, k]=1.$$
Ceci montre que $[Gl_{n}]$ est inversible dans $K(\mathcal{CH}^{sp}(k))$, et 
donc que toutes les classes $\mathbb{L}$ et $\mathbb{L}^{i}-1$ le sont aussi (voir le lemme 
\ref{lfait} ci-dessous).
On consid\`ere alors
l'anneaux localis\'e en $\mathbb{L}$ et en tous les
$\mathbb{L}^{i}-\mathbf{1}$ pour $i>0$, et le
morphisme induit
$$j : K(\mathcal{V}(k))[\mathbb{L}^{-1},\{(\mathbb{L}^{i}-\mathbf{1})^{-1}\}_{i>0}]
\longrightarrow 
K(\mathcal{CH}^{sp}(k)).$$

Notre th\'eor\`eme principal est le suivant. 

\begin{thm}\label{t1}
Le morphisme 
$$j : K(\mathcal{V}(k))[\mathbb{L}^{-1},\{(\mathbb{L}^{i}-\mathbf{1})^{-1}\}_{i>0}]
\longrightarrow 
K(\mathcal{CH}^{sp}(k)).$$
est un isomorphisme. 
\end{thm}

La preuve de ce th\'eor\`eme nous prendra un certain temps et nous
y consacrons le paragraphe suivant. Avant de nous y lancer signalons le fait suivant qui nous sera utile par la
suite. 

\begin{lem}\label{lfait}
On a dans $K(\mathcal{V}(k))$ l'\'egalit\'e suivante
$$[Gl_{n}]=\mathbb{L}^{\frac{n(n-1)}{2}}.\prod_{0<i<n+1} (\mathbb{L}^{i}-\mathbf{1}).$$
En cons\'equence le morphisme naturel 
$$K(\mathcal{CH}^{sp}(k))[\{[Gl_{n}]^{-1}\}_{n>0}] \longrightarrow 
K(\mathcal{CH}^{sp}(k))[\mathbb{L}^{-1},\{(\mathbb{L}^{i}-\mathbf{1})^{-1}\}_{i>0}]$$
est un isomorphisme. 
\end{lem}

\textit{Preuve:} On d\'emontre l'\'egalit\'e
$$[Gl_{n}]=\mathbb{L}^{\frac{n(n-1)}{2}}.\prod_{0<i<n+1}(\mathbb{L}^{i}-\mathbf{1})$$
par r\'ecurrence sur $n$. Pour $n=1$ c'est \'evident. Pour $n>1$, on consid\`ere
l'action naturelle de $Gl_{n}$ sur $\mathbb{A}^{n}-\{0\}$. Cette action est transitive et on dispose donc d'un 
isomorphisme de sch\'emas sur $Spec\, k$
$$\mathbb{A}^{n}-\{0\}\simeq Gl_{n}/H,$$
o\`u $H$ est le stabilisateur de $(1,0,\dots,0)$. Ainsi, le groupe $H$ est le produit semi-direct de
$Gl_{n-1}$ par $\mathbb{G}_{a}^{n-1}$. En particulier, pour tout sch\'ema $X$
le morphisme naturel
$$H^{1}_{zar}(X,H) \longrightarrow H^{1}_{fppf}(X,H)$$
est un isomorphisme (Hilbert 90 et le fait que $\mathbb{G}_{a}^{n-1}$ soit le faisceau en groupes
sous-jacent \`a un faisceau coh\'erent). Le morphisme $Gl_{n} \longrightarrow Gl_{n}/H$ est donc une
$H$-fibration Zariski localement triviale et l'on trouve donc
$$[Gl_{n}]=[H].[Gl_{n}/H]=\mathbb{L}^{n-1}.[Gl_{n-1}].(\mathbb{L}^{n}-\mathbf{1}).$$
On termine alors par r\'ecurrence. 
\hfill $\Box$ \\

\subsection{Preuve du th\'eor\`eme \ref{t1}}

La preuve du th\'eor\`eme \ref{t1} consiste \`a r\'eduire successivement 
la cat\'egorie $\mathcal{CH}^{sp}(k)$, sans changer son groupe de Grothendieck, 
en une chaine qui commence avec $\mathcal{CH}^{sp}(k)$ et qui finit avec
$\mathcal{V}(k)$. \\

Pour commencer, une sous-cat\'egorie pleine 
$\mathcal{C}$ est dite \emph{adapt\'ee} si elle
v\'erifie les conditions suivantes. 

\begin{enumerate}
\item $\mathcal{V}(k) \subset \mathcal{C}$
\item La sous-cat\'egorie
$\mathcal{C}$ est stable par sommes disjointes et produits finis.

\item Si $F\in \mathcal{C}$ alors tout sous-champ ouvert et tout sous-champ  
ferm\'e de $F$ est dans $\mathcal{C}$. 

\item Soit $F_{0}$ un champ qui est soit un sch\'ema affine, soit un 
champ de la forme $K(\mathbb{G}_{a},n)$ pour un entier $n>0$.  Si $F \longrightarrow F'$ est un morphisme
de $\mathcal{C}$ qui est une $F_{0}$-fibration localement Zariski triviale,
alors $F_{0} \in \mathcal{C}$. 
\end{enumerate}

Pour une sous-cat\'egorie pleine
$\mathcal{C}\subset \mathcal{CH}^{sp}(k)$ adapt\'ee, nous noterons
$K(\mathcal{C})$ le groupe quotient 
de $\mathbb{Z}[\mathcal{C}]$ (le groupe ab\'elien libre sur les
classes d'isomorphismes de $\mathcal{C}$) par les trois
relations suivantes.

\begin{enumerate}
\item Pour $F$ et $F'$ dans $\mathcal{C}$ on a
$$[F\coprod F']=[F]+[F'].$$
\item Pour toute $g$-\'equivalence dans $\mathcal{C}$
$F \longrightarrow F'$, on a 
$$[F]=[F'].$$ 
\item Soit $F_{0}$ un champ qui est soit un sch\'ema affine, soit un 
champ de la forme $K(\mathbb{G}_{a},n)$ pour un entier $n>0$. Si $F \longrightarrow F'$ dans $\mathcal{C}$ 
est une $F_{0}$-fibration Zariski localement triviale, alors on a
$$[F]=[F'\times F_{0}].$$
\end{enumerate}

Noter que les relations $(1)$ et $(3)$  ont un sens car 
$\mathcal{C}$ est adapt\'ee. De plus, l'inclusion naturelle
$\mathcal{C} \longrightarrow \mathcal{CH}^{sp}(k)$
induit un morphisme d'anneaux
$$K(\mathcal{C}) \longrightarrow K(\mathcal{CH}^{sp}(k)).$$
Plus g\'en\'eralement, pour deux sous-cat\'egories
adapt\'ees $\mathcal{C}\subset \mathcal{C}'$, le foncteur d'inclusion
induit un morphisme d'anneaux
$$K(\mathcal{C}) \longrightarrow K(\mathcal{C}').$$

Un lemme cl\'e que nous utiliserons tout au long de la preuve est le suivant.

\begin{lem}\label{lcle}
Soit $\mathcal{C} \subset \mathcal{C}'$ deux sous-cat\'egories
pleines et adapat\'ees de $\mathcal{CH}^{sp}(k)$. On suppose que pour 
tout sch\'ema affine $X$, et toute $X$-fibration Zariski localement triviale
$f : F \longrightarrow F'$ dans $\mathcal{C}'$, 
il existe un sous-champ ferm\'e 
$F'_{0}\subsetneq F'$ avec 
$$F'-F'_{0}\in \mathcal{C} \qquad F\times^{h}_{F'}(F'-F_{0}') \in \mathcal{C}.$$ 
Alors le morphisme d'inclusion
$$K(\mathcal{C}) \longrightarrow K(\mathcal{C}')$$
est un isomorphisme. 
\end{lem}

\textit{Preuve:} On d\'efinit un inverse 
$$\phi : K(\mathcal{C}') \longrightarrow K(\mathcal{C})$$
de la fa\c{c}on suivante. Pour $F \in \mathcal{C}'$, on peut trouver une
suite d\'ecroissante de sous-champs ferm\'es 
$$\xymatrix{\emptyset=F_{r+1} \subset F_{r} \subset \dots \subset 
F_{1}\subset F_{0}=F}$$
tel que chaque champ $F_{i}-F_{i+1}$ soit dans $\mathcal{C}$ (on applique l'hypoth\`ese
\`a $X=Spec\, k$). 
On pose 
$$\phi([F]):=\Sigma_{i}[F_{i}-F_{i+1}].$$
Il faut montrer que ceci est bien d\'efini, c'est \`a dire compatible avec les relations
qui d\'efinissent $K(\mathcal{C}')$ et  $K(\mathcal{C}')$.  Pour les deux premi\`eres 
relations ceci est imm\'ediat. La troisi\`eme relation est compatible grace \`a l'hypoth\`ese
(on remarquera que pour une $K(\mathbb{G}_{a},n)$-fibration Zariski localement triviale
$F \longrightarrow F'$
avec $n>0$, les sous-champs ouverts de $F'$ sont en correspondance avec ceux
de $F$). 
Enfin, on v\'erifie facilement que $\phi$ est l'inverse cherch\'e. 
\hfill $\Box$\\

Soit $\mathcal{C}_{1}$ la sous-cat\'egorie pleine
de $\mathcal{CH}^{sp}(k)$ form\'ee des champs 
de $F$ tel qu'il existe un espace alg\'ebrique $X$, un entier
$r\geq 0$ et un 
$Gl_{r}$-torseur $X \longrightarrow \tau_{\leq 1}F$ (on posera
$Gl_{0}:=Spec\, k$). C'est une
sous-cat\'egorie adapt\'ee. Notons que si $F\in \mathcal{C}_{1}$ alors 
le $1$-champ $\tau_{\leq 1}F$ est \'equivalent au champ quotient $[X/Gl_{r}]$, et donc est un $1$-champ d'Artin. 
Noter aussi que d'apr\`es le th\'eor\`eme de Hilbert 90, tout $Gl_{r}$-torseur entre champs est une
$Gl_{r}$-fibration Zariski localement triviale. 
 
\begin{lem}\label{l3}
Le morphisme
$$K(\mathcal{C}_{1}) \longrightarrow K(\mathcal{CH}^{sp}(k))$$
est un isomorphisme. 
\end{lem}

\textit{Preuve:} Pour pouvoir appliquer le lemme \ref{lcle} il nous
suffit de montrer que pour $F\in \mathcal{C}_{1}$, il existe un
sous-champ ferm\'e $F_{0}$ de $F$ tel que $F-F_{0} \in \mathcal{C}_{1}$.
En effet, supposons que ceci soit le cas. Alors, si $f : F \longrightarrow F'$ est une $X$-fibration localement Zariski triviale avec 
$X$ un sch\'ema affine, on v\'erifie facilement que le morphisme induit
$\tau_{\leq 1}F \longrightarrow \tau_{\leq 1}F'$ est encore une
$X$-fibration localement Zariski triviale. De ceci on d\'eduit que
$F$ est dans $\mathcal{C}_{1}$ si $F'$ l'est, et donc que les conditions
du lemme \ref{lcle} sont v\'erifi\'ees. 

Montrons donc que pour $F\in \mathcal{C}_{1}$, il existe un
sous-champ ferm\'e $F_{0}$ de $F$ tel que $F-F_{0} \in \mathcal{C}_{1}$.
En rempla\c{c}ant $F$ par son sous-champ r\'eduit on voit que l'on 
peut supposer que $F$ est r\'eduit. Il s'agit alors de trouver
un ouvert non-vide $U$ de $F$ tel que $U \in \mathcal{C}_{1}$. 
Par le corollaire \ref{c1'} on voit aussi que l'on peut supposer que $F$ est une gerbe
totale, et en particulier que $\tau_{\leq 1}F$ est un $1$-champ d'Artin. 
Comme ce que l'on cherche \`a montrer ne concerne que le tronqu\'e $\tau_{\leq 1}F$ on peut
tout simplement supposer que $F=\tau_{\leq 1}F$ (i.e. que
$F$ est 1-tronqu\'e).

On consid\`ere la projection
$$F \longrightarrow M(F)=M.$$
On restreint ce morphisme \`a un point g\'en\'erique $Spec\, K$ de $M$, et on 
trouve une gerbe
$F_{0}$ au-dessus de $Spec\, K$. Il existe une extension finie
$K'$ de $K$, tel que le champ $F_{0}\times^{h}_{Spec\, K}Spec\, K'$ soit de le forme
$K(G',1)$, o\`u $G'$ est un sch\'ema en groupes affine et de type fini sur $K'$. 
On choisit alors une immersion ferm\'ee de sch\'emas en groupes sur $K'$
$$i : G' \hookrightarrow Gl_{r'}.$$
Le morphisme $i$ d\'efinit un fibr\'e vectoriel $E$ de rang $r'$ sur $K(G',1)$. Le morphisme
$$p : K(G',1) \longrightarrow F_{0}$$
\'etant fini et plat, $f_{*}(E)$ est un fibr\'e vectoriel sur $F_{0}$, de rang disons $r$. Le $Gl_{r}$-torseur associ\'e 
d\'efini un morphisme de champs $X_{0} \longrightarrow F_{0}$. Par construction on peut voir que
$X_{0}$ est $0$-tronqu\'e (car $i$ est un monomorphisme), et donc que c'est un espace
alg\'ebrique.  

Maintenant, comme $K(Gl_{r},1)$ est un champ localement de pr\'esentation finie, le $Gl_{r}$-torseur
$X \longrightarrow F_{0}$ s'\'etend en un  
$Gl_{r}$-torseur $X \longrightarrow F\times^{h}_{M}U$, pour un ouvert non-vide
$U$ de $M$. En choisissant $U$ assez petit on peut aussi faire en sorte que $X$ soit
un espace alg\'ebrique.  \hfill $\Box$ \\

Notons $\mathcal{C}_{2}$ la sous-cat\'egorie pleine de 
$\mathcal{C}_{1}$ form\'ee des champs $F$ tels que 
pour tout  sch\'ema affine $X$ et tout point $x \in F(K)$, 
les faisceaux $\pi_{i}(F,x)$ sur $X$ soient 
repr\'esentables par des sch\'emas en groupes plats sur $X$. 
C'est une sous-cat\'egorie adapt\'ee. 

\begin{lem}\label{l4}
Le morphisme
$$K(\mathcal{C}_{2}) \longrightarrow 
K(\mathcal{C}_{1})$$
est un isomorphisme. 
\end{lem}

\textit{Preuve:} C'est une application du lemme cl\'e \ref{lcle} et du corollaire
\ref{c1'}. \hfill $\Box$ \\

Soit $\mathcal{C}_{3}$ la sous-cat\'egorie pleine de 
$\mathcal{C}_{2}$ form\'ee des champs $F$ tels que le morphisme naturel
$F\longrightarrow \tau_{\leq 1}F$ soit un isomorphisme (i.e. form\'ee des
$1$-champs d'Artin).  C'est une sous-cat\'egorie adapt\'ee. 

\begin{lem}\label{l5}
Le morphisme d'inclusion 
$$K(\mathcal{C}_{3}) \longrightarrow 
K(\mathcal{C}_{2})$$
est injectif. 
\end{lem}

\textit{Preuve:} Nous allons construire une morphisme
$$\phi : K(\mathcal{C}_{2})
\longrightarrow K(\mathcal{C}_{3})$$
et montrer que c'est une retraction. 

Soit $F\in \mathcal{C}_{2}$, et consid\'erons le morphisme
naturel $F \longrightarrow \tau_{\leq 1}F$. Ce morphisme
est relativement simplement connexe, et donc on peut d\'efinir
des faisceaux en groupes ab\'eliens $\pi_{i}(F)$ sur
$\tau_{\leq 1}F$. Ces faisceaux sont tels que pour tout
sch\'ema affine $X$ et tout morphisme $x : X \longrightarrow F$, 
la restriction de $\pi_{i}(F)$ sur $X$ soit naturellement isomorphe 
au faisceau $\pi_{i}(F,x)$. En particulier on voit
que les $\pi_{i}(F)$ sont des sch\'emas en groupes plats et unipotents sur
$\tau_{\leq 1}(F)$. Ils poss\`edent donc une dimension relative
localement constante, et on peut donc \'ecrire
$\tau_{\leq 1}F$ comme une r\'eunion disjointe de 
sous-champs ouverts et ferm\'es $\tau_{\leq 1}F_{\alpha} \subset \tau_{\leq 1}F$ tel que chaque sch\'ema
en groupes $\pi_{i}(F)$ sur $\tau_{\leq 1}F_{\alpha}$ soit 
de dimension relative constante \'egale \`a 
$\beta_{i}^{F,\alpha}$ (par la suite on verra
$\beta_{i}^{F,\alpha}$ comme la valeur d'une fonction 
$\beta_{i}^{F}$ localement constante sur $\tau_{\leq 1}F$). On pose alors
$$\phi([F]):=\sum_{\alpha}[\tau_{\leq 1}F_{\alpha}].\prod_{i>1}.\mathbb{L}^{(-1)^{i}.\beta_{i}^{F,\alpha}}
\in K(\mathcal{C}_{3}).$$

Il est imm\'ediat de v\'erifier que $\phi$ ainsi d\'efini est compatible avec la
relation $(1)$ de la d\'efinition \ref{d6}. Soit $f : F \longrightarrow F'$ 
une g-\'equivalence. On voit que pour tout corps alg\'ebriquement 
clos $K$ et tout point $x \in F(K)$, le morphisme induit
$$\pi_{i}(F,x)(K)\simeq \pi_{i}(F(K),x) \longrightarrow \pi_{i}(F',f(x))(K)\simeq \pi_{i}(F'(K),f(x))$$
est un isomorphisme. Ceci implique clairement que 
les sch\'emas en groupes $\pi_{i}(F,x)$ et $\pi_{i}(F',f(x))$ ont m\^eme dimension. On peut ainsi \'ecrire
$F'$ comme r\'eunion disjoint $\coprod_{\alpha}F'_{\alpha}$, et de m\^eme 
$F\simeq \coprod_{\alpha}F_{\alpha}$ avec $F_{\alpha}\simeq F\times_{F'}^{h}F'_{\alpha}$, 
tel que $\beta_{i}^{F,\alpha}=\beta_{i}^{F',\alpha}$ pour tout $i>1$. Comme de plus
le morphisme induit $\tau_{\leq 1}F \longrightarrow \tau_{\leq 1}F'$ est une g-\'equivalence, ceci implique clairement que
$\phi([F])=\phi([F'])$. On voit donc que $\phi$ est aussi compatible avec la relation $(2)$ 
de la d\'efinition \ref{d6}. 

Montrons enfin que $\phi$ est aussi compatible avec 
la relation $(3)$. Soit $F_{0}$ un champ qui est soit $K(\mathbb{G}_{a},n)$ pour
un $n>0$, soit un un sch\'ema affine, et soit $f : F \longrightarrow F'$ une $F_{0}$-fibration 
Zariski localement triviale dans $\mathcal{C}_{2}$. Dans le cas o\`u 
$F_{0}$ est un sch\'ema affine on voit facilement que 
l'on a une \'equalit\'e $f^{*}(\beta_{i}^{F'})=\beta_{i}^{F}$ pour tout $i>1$. Comme de plus
le morphisme induit 
$\tau_{\leq 1}F \longrightarrow \tau_{\leq 1}F'$ est une $F_{0}$-fibration Zariski localement triviale
ceci implique que 
$\phi([F])=\phi([F_{0}]).\phi([F'])$. 

Supposons maintenant que 
$F_{0}$ soit de la forme $K(\mathbb{G}_{a},n)$ pour $n>1$, ce qui
implique que $\tau_{\leq 1}F\simeq \tau_{\leq 1}F'$.  
En utilisant la relation $(1)$ de la d\'efinition \ref{d6} on peut supposer que
$F$ et $F'$ sont tous deux connexes. On note alors simplement $\beta_{i}$ la valeur 
de $\beta_{i}^{F}$ et $\beta_{i}'$ la valeur de $\beta_{i}^{F'}$. La suite exacte longue en homotopie
montre que pour $i\neq n,n+1$, le morphisme induit $\pi_{i}(F) \longrightarrow \pi_{i}(F')$
est un isomorphisme. Ceci implique que pour $i\neq n,n+1$ on a
$\beta_{i}=\beta_{i}'$. De plus, on a une suite exacte longue de faisceaux en groupes ab\'eliens sur $\tau_{\leq 1}F\simeq
\tau_{\leq 1}F'$
$$\xymatrix{
0 \ar[r] & \pi_{n+1}(F) \ar[r] &  \pi_{n+1}(F') \ar[r] & \mathbb{G}_{a} \ar[r]^-{u} & \pi_{n}(F) \ar[r] &
 \pi_{n}(F') \ar[r] & 0.}$$
Rappelons que la dimension des fibres des $\pi_{i}(F)$ et des $\pi_{i}(F')$
est constante. Ainsi, suivant que le morphisme $u$ est nul ou pas, on voit que l'on a deux
cas, ou bien 
$$\beta_{n+1}'=\beta_{n+1}+1 \qquad \beta_{n}'=\beta_{n},$$
ou bien 
$$\beta_{n+1}'=\beta_{n+1} \qquad \beta_{n}=\beta_{n}'+1.$$
Dans tous les cas ceci implique que l'on a 
$$\phi([F])=\phi([F']).\mathbb{L}^{(-1)^{n}}=\phi([F' \times K(\mathbb{G}_{a},n)]).$$
Il nous reste le cas o\`u $F_{0}=K(\mathbb{G}_{a},1)$. On dispose alors d'une suite exacte 
de faisceaux sur $\tau_{\leq 1}F$ 
$$\xymatrix{0 \ar[r] & \pi_{2}(F) \ar[r] & f^{*}(\pi_{2}(F')) \ar[r]^-{\delta} & \mathbb{G}_{a}\times \tau_{\leq 1}F.}$$ 
Ainsi, comme les dimensions relatives de $\pi_{2}(F)$ et 
$\pi_{2}(F')$ sont constantes on voit que ou bien 
le morphisme $\delta$ est surjectif, ou bien l'image de $\delta$ est fibres \`a fibres
un sous-groupe de dimension nulle de $\mathbb{G}_{a}$. Dans le premier cas, on a
$$\beta_{2}=\beta_{2}'-1 \qquad \beta_{i}=\beta_{i}' \; \;  \forall \; i\neq 2$$
et $\tau_{\leq 1}F\simeq \tau_{\leq 1}F'$. Ceci implique clairement que
$\phi([F])=\phi([F_{0}]).\phi([F'])$. Dans le second
cas on a $\beta_{i}=\beta_{i}'$ pour tout $i>1$, et donc il nous suffit de voir que
$[\tau_{\leq 1}F]=[\tau_{\leq 1}F'].\mathbb{L}^{(-1)}$ dans 
$K(\mathcal{C}_{3})[\mathbb{L}^{-1}]$. Pour cela notons que 
$\tau_{\leq 1}F \longrightarrow \tau_{\leq 1}F'$ est un morphisme de $1$-champs d'Artin tel que 
pour tout sch\'ema affine $X$ et tout point $x : X \longrightarrow F'$, le champ
$F\times_{F'}^{h}X$ est, localement sur $X$, \'equivalent comme champ sur $X$ \`a un champ de la forme 
$K(H,1)$, avec $H$ un sch\'ema en groupes sur $X$, 
\`a fibres unipotentes connexes et de dimension
$1$. Ceci implique que, tout au moins sur un ouvert non vide de $\tau_{\leq 1}(F')$ (ce qui suffit
par un argument de stratification), 
que le morphisme $\tau_{\leq 1}F \longrightarrow \tau_{\leq 1}F'$ est une $K(\mathbb{G}_{a},1)$-fibration 
Zariski localement triviale. Ainsi, on a 
$[\tau_{\leq 1}F]=[\tau_{\leq 1}F'].[K(\mathbb{G}_{a},1)]=[\tau_{\leq 1}F'].\mathbb{L}^{(-1)}$. 
Ceci termine la d\'emonstration du fait que
$\phi$ soit compatible avec la relation $(3)$ de la d\'efinition \ref{d6}, et donc qu'il d\'efinisse un morphisme
$$K(\mathcal{C}_{2})
\longrightarrow K(\mathcal{C}_{3}).$$
Par construction $\phi$ est une r\'etraction du morphisme d'inclusion 
$$i : K(\mathcal{C}_{3}) \longrightarrow 
K(\mathcal{C}_{2}).$$
\hfill $\Box$ \\

\begin{lem}\label{l6}
Le morphisme d'inclusion 
$$i : K(\mathcal{C}_{3})
\longrightarrow K(\mathcal{C}_{2})$$
est surjectif. 
\end{lem}

\textit{Preuve:}
Soit donc $F\in \mathcal{C}_{2}$ (que l'on peut supposer r\'eduit). Soit 
$X \longrightarrow \tau_{\leq 1}F$ un $Gl_{r}$-torseur avec $X$ un espace alg\'ebrique. On pose 
$F_{X}:=F\times^{h}_{\tau_{\leq 1}F}X$. Alors, le morphisme $F_{X} \longrightarrow F$ 
est un $Gl_{r}$-torseur, et donc une $Gl_{r}$-fibration Zariski localement triviale. Ainsi
on a $[F]=[F_{X}].[Gl_{r}]^{-1}$ dans $K(\mathcal{C}_{2})$.
On voit ainsi que l'on peut remplacer $F$ par $F_{X}$, o\`u encore ce qui revient au m\^eme supposer que
$F$ est une gerbe totale avec $\tau_{\leq 1}F\simeq M(F)=:M$. 

On suppose maintenant que $F$ est un $n$-champ
d'Artin, et on montre par r\'ecurrence sur $n$ que $[F]$ est dans l'image 
du morphisme $i$. Pour $n=0,1$ c'est \'evident car $F\in \mathcal{C}_{3}$. 
Supposons donc $n>1$, et que pour tout $(n-1)$-champ $F \in \mathcal{C}_{2}$,
$[F]$ soit dans l'image de $i$. On consid\`ere le morphisme naturel
$$p : F \longrightarrow \tau_{\leq n-1}F.$$
C'est un fibration localement triviale pour la topologie fppf de fibre
\'equivalente \`a $K(H,n)$, o\`u $H=\pi_{n}(F)$ est un sch\'ema en groupes
ab\'eliens plats et unipotent sur $M(F)=M$. Sur un ouvert Zariski
$U$ non vide de $M$, il existe une suite exacte courte 
de sch\'emas en groupes sur $M$
$$\xymatrix{
0 \ar[r] & H_{red} \ar[r] & H \ar[r] & K \ar[r] & 0,}$$
o\`u $H_{red}$ est r\'eduit, et $K:=H/H_{red}$ est plat et de dimension relative nulle sur $M$. En stratifiant $F$ 
si n\'ecessaire
on supposera donc qu'une telle suite exacte existe sur $M$. On pourra aussi, et c'est ce que nous
ferons, supposer que le groupe $H_{red}$ poss\`ede une fitration 
$$H_{0} \subset H_{1} \subset \dots \subset H_{r}=H_{red},$$
par des sous-groupes plats sur $M$ et telle que chaque
quotient $H_{i}/H_{i-1}$ soit isomorphe, comme sch\'ema en groupes sur $M$, \`a $\mathbb{G}_{a}\times M$.

Maintenant, la fibration $p$ est 
classifi\'ee par une classe dans $H^{n+1}_{fppf}( \tau_{\leq n-1}F,H)$, et entre donc dans
un diagramme homotopiquement cart\'esien
$$\xymatrix{
F \ar[r] \ar[d] &  \tau_{\leq n-1}F \ar[d] \\
M \ar[r] & K(H,n+1).}$$
En composant avec le morphisme naturel $H \longrightarrow K$, on trouve 
un autre diagramme homotopiquement cart\'esien
$$\xymatrix{
F' \ar[r] \ar[d] &  \tau_{\leq n-1}F \ar[d] \\
M \ar[r] & K(K,n+1).}$$
Il existe de plus un morphisme naturel $F \longrightarrow F'$, qui  est 
une fibration en $K(H_{red},n)$, et il existe donc un diagramme homotopiquement cart\'esien
$$\xymatrix{
F \ar[r] \ar[d] &  F' \ar[d] \\
M \ar[r] & K(H_{red},n+1).}$$
En composant avec les quotients successifs $H \longrightarrow H/H_{i}$, on d\'efinit
des champs $F_{i}$  par
des diagrammes homotopiquement cart\'esiens
$$\xymatrix{
F_{i} \ar[r] \ar[d] &  F' \ar[d] \\
M \ar[r] & K(H/H_{i},n+1).}$$
De plus, il existe des morphismes naturels $F_{i-1} \longrightarrow F_{i}$ qui sont 
des fibrations localement triviales en 
$K(\mathbb{G}_{a},n)$. 
Ainsi, en utilisant la relation
$(3)$ de la d\'efinition \ref{d6} on trouve dans 
$K(\mathcal{C}_{2})$
$$[F_{i-1}]=[F_{i}].[K(\mathbb{G}_{a},n)]=[F_{i}].\mathbb{L}^{(-1)^{n}}.$$
Par r\'ecurrence sur $i$ ceci donne 
$$[F]=[F'].\mathbb{L}^{(-1)^{n}.d},$$
o\`u $d$ est la dimension relative de $H$ sur $M$. 
Ainsi, on peut supposer que 
$F=F'$, et donc que $H$ est un groupe fini unipotent et plat sur $M$.

Mais pour tout corps $K$ alg\'ebriquement 
clos on a pour tout $i>0$
$$\pi_{i}(K(H,n)(K))\simeq H^{n-i}_{fppf}(Spec\, K,H)\simeq 0.$$
Ainsi, l'ensemble simplicial $K(H,n)(K)$ est contractile, et donc
le morphisme $F\longrightarrow \tau_{\leq n-1}F$
est une g-\'equivalence. On a donc
$[F]=[\tau_{\leq n-1}F]$ dans $K(\mathcal{C}_{2})$, ce qui 
termine la preuve du lemme.
\hfill $\Box$ \\

Soit $\mathcal{C}_{4}$ la sous-cat\'egorie pleine 
de $\mathcal{C}_{3}$ form\'ee des espaces alg\'ebriques

\begin{lem}\label{l7}
Le morphisme d'inclusion
$$K(\mathcal{C}_{4})[\mathbb{L}^{-1},\{(\mathbb{L}^{i}-\mathbf{1})^{-1}\}_{i>0}] \longrightarrow 
K(\mathcal{C}_{3})$$
est un isomorphisme. 
\end{lem}

\textit{Preuve:}  Rappelons que $\mathcal{C}_{3}$, qui par d\'efinition est une sous-cat\'egorie
de $\mathcal{C}_{1}$, est la cat\'egorie des $1$-champs d'Artin $F$ tel qu'il existe une espace
alg\'ebrique $X$ et un $Gl_{r}$-torseur $X \longrightarrow F$, ou de mani\`ere \'equivalente tel que
$F$ soit \'equivalent \`a un champ quotient $[X/Gl_{r}]$ pour un certain espace alg\'ebrique
$X$ et une certaine action de $Gl_{r}$ sur $X$. On construit alors un morphisme 
$$\phi : K(\mathcal{C}_{3}) \longrightarrow 
K(\mathcal{C}_{4})[\mathbb{L}^{-1},\{(\mathbb{L}^{i}-\mathbf{1})^{-1}\}_{i>0}]$$
en posant $\phi([F])=[X].[Gl_{r}]^{-1}$ (cette formule a un sens grace au lemme \ref{lfait}). 
On remarque facilement que cette d\'efinition est ind\'ependant du choix
du $Gl_{r}$-torseur $X \longrightarrow F$ (\`a l'aide de la relation $(3)$ 
de la d\'efinition \ref{d6}). 

Il nous faut montrer que $\phi$ ainsi d\'efini est compatible aux trois
relations de la d\'efinition \ref{d6}. La premi\`ere est \'evidente. Soit 
$f : F \longrightarrow F'$ une $g$-\'equivalence dans 
$\mathcal{C}_{3}$. On choisit un $Gl_{r}$-torseur $X' \longrightarrow F'$ avec $X'$ un espace alg\'ebrique, et on
consid\`ere le carr\'e homotopiquement cart\'esien suivant
$$\xymatrix{
F \ar[r]  & F'  \\
Y \ar[r] \ar[u] & X'\ar[u].
}$$
Le morphisme $Y \longrightarrow F$ est toujours un $Gl_{r}$-torseur, mais $Y$ n'est peut-\^etre pas
un espace alg\'ebrique. On choisit donc un $Gl_{s}$-torseur $X \longrightarrow F$ avec 
$X$ un espace alg\'ebrique, et on consid\`ere
$Z:=X\times_{F}^{h}Y$. Le morphisme naturel $Z \longrightarrow X$ est un $Gl_{r}$-torseur, et
donc $Z$ est un espace alg\'ebrique. De plus le morphisme $Z \longrightarrow Y$ est un 
$Gl_{s}$-torseur. Enfin, comme le morphisme $Y \longrightarrow X'$ est une $g$-\'equivalence, on trouve dans
$K(\mathcal{C}_{4})[\mathbb{L}^{-1},\{(\mathbb{L}^{i}-\mathbf{1})^{-1}\}_{i>0}]$ l'\'egalit\'e suivante
$$[X'].[Gl_{r}]^{-1}=[Y].[Gl_{r}]^{-1}=[Z].[Gl_{s}]^{-1}.[Gl_{r}]^{-1}=[X].[Gl_{s}]^{-1}.$$
Ceci implique bien que $\phi([F])=\phi([F'])$. Il nous reste donc \`a v\'erifier que 
$\phi$ est compatible avec la relation $(3)$ de la d\'efinition \ref{d6}. 
Pour cela, soit $F \longrightarrow F'$ une $F_{0}$-fibration Zariski localement triviale, avec
$F_{0}$ qui est soit un sch\'ema affine $X$ soit $K(\mathbb{G}_{a},1)$ (noter que comme
$F$ et $F'$ sont des $1$-champs on ne peut pas avoir une fibration 
en $K(\mathbb{G}_{a},n)$ pour $n>1$). Commen\c{c}ons par le cas o\`u $F_{0}$ est un sch\'ema affine
$Z$. Soit $X' \longrightarrow F'$ un $Gl_{r}$-torseur avec $X'$ un espace alg\'ebrique, et $X:=X'\times_{F'}^{h}F \longrightarrow F$ 
le $Gl_{r}$-torseur induit. Comme le morphisme $X \longrightarrow X'$ est une $Z$-fibration 
Zariski localement triviale on voit que $X$ est un espace alg\'ebrique. Ainsi, on trouve
$$[X].[Gl_{r}]^{-1}=[Z].[X'].[Gl_{r}]^{-1},$$
ce qui montre que $\phi([F])=\phi([F']).\phi([Z]))$. Passons au cas o\`u $F_{0}$ est 
$K(\mathbb{G}_{a},1)$. Notons encore $X' \longrightarrow F'$ un $Gl_{r}$-torseur avec $X'$ un espace
alg\'ebrique et $Y:=X'\times_{F'}^{h}F$. Le morphisme $Y \longrightarrow F$ est un $Gl_{r}$-torseur, 
et $Y \longrightarrow X'$ est une $K(\mathbb{G}_{a},1)$-fibration Zariski localement triviale. 
Soit $X \longrightarrow F$ un $Gl_{s}$-torseur, et $Z:=X\times_{F}^{h}Y$. Alors, 
comme pr\'ec\'edemment on voit que $Z$ est un espace alg\'ebrique et que
$Z \longrightarrow Y$ est un $Gl_{s}$-torseur. Soit $X'_{i}\subset X'$ une suite d\'ecroissante
de sous-espaces alg\'ebrique ferm\'es de $X'$ tel que pour $V_{i}=X'_{i}-X'_{i-1}$ le
champ $Y_{i}:=Y\times^{h}_{X'}V_{i}$ soit \'equivalent au-dessus de $V_{i}$ \`a
la projection $V_{i}\times K(\mathbb{G}_{a},1)$. Posons alors
$Z_{i}:=Z\times_{Y}^{h}Y_{i}$ et $W_{i}:=Z\times_{Y}^{h}V_{i}$. Alors
les morphismes  $W_{i} \longrightarrow V_{i}$ sont des $Gl_{s}$-torseurs entre espaces alg\'ebriques. 
D'autre part les morphismes  $W_{i} \longrightarrow Z_{i}$
sont des $\mathbb{G}_{a}$-torseurs entre espaces alg\'ebrique. Ainsi, on trouve les \'egalit\'es suivantes
dans $K(\mathcal{C}_{4})[\mathbb{L}^{-1},\{(\mathbb{L}^{i}-\mathbf{1})^{-1}\}_{i>0}]$

$$[X].[Gl_{s}]^{-1}=[Z].[Gl_{r}]^{-1}.[Gl_{s}]^{-1}=
\Sigma_{i}  [Z_{i}].[Gl_{r}]^{-1}.[Gl_{s}]^{-1}=\Sigma_{i}  [W_{i}].\mathbb{L}^{-1}.[Gl_{r}]^{-1}.[Gl_{s}]^{-1}$$
$$=\Sigma_{i}  [V_{i}].\mathbb{L}^{-1}.[Gl_{r}]^{-1}=[X'].[Gl_{r}]^{-1}.\mathbb{L}^{-1}.$$
En d'autres termes on trouve
$$\phi([F])=\phi([F']).\mathbb{L}^{-1}.$$
Ceci est presque ce que l'on cherche \`a d\'emontrer, pour finir il nous reste encore \`a remarquer que
$\phi([K(\mathbb{G}_{a},1)])=\mathbb{L}^{-1}$. Pour cela, soit 
$\mathbb{G}_{a} \subset Gl_{2}$ le sous-groupe des matrices triangulaires sup\'erieures
avec des $1$ sur la diagonale. Si l'on pose $X=Gl_{2}/\mathbb{G}_{a}$, le morphisme naturel
$X \longrightarrow K(\mathbb{G}_{a},1)$ est un $Gl_{2}$-torseur. On trouve donc
$\phi([K(\mathbb{G}_{a},1)])=[X].[Gl_{2}]^{-1}$.
Cependant, le morphisme $Gl_{2} \longrightarrow X$ est un $\mathbb{G}_{a}$-torseur, et donc
une $\mathbb{G}_{a}$-fibration Zariski localement triviale. On a donc bien
$\phi([K(\mathbb{G}_{a},1)])=[X].[Gl_{2}]^{-1}=[\mathbb{G}_{a}]^{-1}.[Gl_{2}].[Gl_{2}]^{-1}=\mathbb{L}^{-1}$.

Nous en avons fini avec le fait que $\phi$ soit compatible aux relations
$(1)$, $(2)$ et $(3)$ de la d\'efinition \ref{d6}, et donc avec le fait que
$\phi$ soit bien d\'efini. De plus,
il est clair que le morphisme $\phi$ est une r\'etraction du morphisme d'inclusion
$$K(\mathcal{C}_{4})[\mathbb{L}^{-1},\{(\mathbb{L}^{i}-\mathbf{1})^{-1}\}_{i>0}] \longrightarrow 
K(\mathcal{C}_{3}).$$
Enfin, pour tout $F \in \mathcal{C}_{3}$, soit $X$ un espace alg\'ebrique et $X \longrightarrow F$ 
un $Gl_{r}$-torseur. Alors, on a dans $K(\mathcal{C}_{3})[\mathbb{L}^{-1},\{(\mathbb{L}^{i}-\mathbf{1})^{-1}\}_{i>0}]$ 
l'\'egalit\'e $[F]=[X].[Gl_{r}]^{-1}$. Ceci montre que le morphisme d'inclusion
est surjectif et donc que c'est un isomorphisme.  \hfill $\Box$ \\

Notons enfin $\mathcal{C}_{5}=\mathcal{V}(k)$ la sous-cat\'egorie pleine form\'ee
des sch\'emas de type fini sur $k$. Le lemme \ref{lcle} implique que le morphisme d'inclusion
$$K(\mathcal{V}(k)) \longrightarrow K(\mathcal{C}_{4})$$
est un isomorphisme. Ainsi, les lemmes \ref{l3}, \ref{l4}, \ref{l5}, \ref{l6} et \ref{l7} mis bout \`a bout montrent que le morphisme
d'inclusion
$$K(\mathcal{V}(k))[\mathbb{L}^{-1},\{(\mathbb{L}^{i}-\mathbf{1})^{-1}\}_{i>0}]
\longrightarrow 
K(\mathcal{CH}^{sp}(k))$$
est un isomorphisme. Ceci ach\`eve la preuve du th\'eor\`eme \ref{t1}. 

\subsection{Invariants num\'eriques des champs d'Artin sp\'eciaux}

Une autre fa\c{c}on d'\'enoncer le th\'eor\`eme \ref{t1} est la suivante. 

\begin{cor}\label{c3}
Soit $A$ un anneau commutatif et 
$$\chi : Iso(\mathcal{V}(k)) \longrightarrow A$$
une application des classes d'isomorphismes de vari\'et\'es sur $k$ vers $A$ qui v\'erfie les  
conditions suivantes. 
\begin{enumerate}
\item $$\chi(X\coprod Y)=\chi(X)+\chi(Y)$$
\item $$\chi(Spec\, k)=1$$
\item $$\chi(X\times Y)=\chi(X).\chi(Y)$$
\item Pour toute g-\'equivalence $g : X \longrightarrow Y$ on 
a $\chi(X)=\chi(Y)$.
\item Les \'el\'ements $\chi(\mathbb{A}^{1})$ et $\chi(\mathbb{A}^{i}-\{0\})$ pour $i>0$ sont 
inversibles dans $A$. 
\end{enumerate}
Alors, il existe une extension de $\chi$ en une application
$$\chi : Iso(\mathcal{CH}^{sp}(k)) \longrightarrow A.$$
Cette extension est de plus unique si elle satisfait les conditions suivantes.
\begin{enumerate}
\item $$\chi(F\coprod F')=\chi(F)+\chi(F')$$
\item $$\chi(F\times F')=\chi(F).\chi(F')$$
\item Pour toute g-\'equivalence $g : F \longrightarrow F'$ on 
a $\chi(F)=\chi(F')$.
\item Pour tout $F_{0}$ qui est soit un sch\'ema affine, soit 
de la forme $K(\mathbb{G}_{a},n)$ avec $n>1$, et toute
$F_{0}$-fibration Zariski localement triviale, on a
$$\chi(F)=\chi(F_{0}).\chi(F').$$
\end{enumerate}
\end{cor}

Le corollaire pr\'ec\'edent permet de construire de tr\`es exemples
de morphismes d'anneaux $\chi : K(\mathcal{CH}^{sp}(k)) \longrightarrow A$, et 
donc de tr\`es nombreux invariants de champs sp\'eciaux. En effet, on connait beaucoup 
d'exemples de morphismes d'anneaux
$$\chi : K(\mathcal{V}(k)) \longrightarrow A$$
qui rendent $\mathbb{L}$ et les $\mathbb{L}^{i}-\mathbf{1}$ inversibles. Citons les exemples suivants. \\

\textit{S\'eries et nombres de Hodge:} Supposons que $k=\mathbb{C}$. Alors il existe
un morphisme d'anneaux
$$P_{H} : K(\mathcal{V}(k)) \longrightarrow 
\mathbb{Z}[u,v],$$
d\'etermin\'e par la propri\'et\'e que pour $X$ une vari\'ete lisse et projective sur $k$ on ait
$$P_{H}(X)(u,v)=\Sigma_{p,q}(-1)^{p+q}Dim H^{p}(X,\Omega^{q}_{X})u^{p}v^{q}$$
(voir \cite{dl}). Noter que l'on a 
$$P_{H}(\mathbb{A}^{1})=P_{H}(\mathbb{P}^{1})-P_{H}(Spec\, k)=1+uv-1=uv.$$
De m\^eme, on a pour $i>0$
$$P_{H}(\mathbb{L}^{i}-\mathbf{1})=(uv)^{i}-1.$$
Ainsi, le th\'eor\`eme \ref{t1} nous dit que le morphisme
$$P_{H} : K(\mathcal{V}(k)) \longrightarrow 
\mathbb{Z}[[u,v]][u^{-1},v^{-1}]$$
qui se factorise en un morphisme
$$P_{H} : K(\mathcal{V}(k))[\mathbb{L}^{-1},\{(\mathbb{L}^{i}-\mathbf{1})^{-1}\}_{i>0}] \longrightarrow 
\mathbb{Z}[[u,v]][u^{-1},v^{-1}]$$
s'\'etend de fa\c{c}on unique en un morphisme d'anneaux
$$P_{H} : K(\mathcal{CH}^{sp}(k)) \longrightarrow 
\mathbb{Z}[[u,v]][u^{-1},v^{-1}].$$

\begin{df}\label{d10}
Pour un champ d'Artin sp\'ecial $F$, \emph{la s\'erie de Hodge de $F$} est
la s\'erie de Laurent en deux variables 
$$P_{H}([F])(u,v)=\Sigma_{p,q}h^{p,q}(F)u^{p}v^{q}\in \mathbb{Z}[[u,v]][u^{-1},v^{-1}].$$
Les nombres entiers $h^{p,q}(F)$ sont appel\'es les
\emph{nombres de Hodge de $F$}.
\end{df}

\textit{Caract\'eristique d'Euler motivique:} Supposons maintenant que 
$k$ soit de caract\'eristique nulle. Notons $CMot(k)$ la cat\'egorie des 
motifs de Chow sur $k$.
 En utilisant la r\'esolution des singularit\'es on peut
consrtuire un morphisme d'anneaux
$$\chi_{mot} : K(\mathcal{V}(k)) \longrightarrow 
K(CMot(k)),$$
o\`u $K(CMot(k))$ est le groupe de Grothendieck de la cat\'egorie 
additive $CMot(k)$ (voir \cite{gs}). Notons que $\chi_{mot}(\mathbb{L})=[L]$, avec 
$L=h^{2}(\mathbb{P}^{1})$ le motif de Lefschetz, et est donc un objet 
inversible de la cat\'egorie des $CMot(k)$. 
Ainsi, le th\'eor\`eme \ref{t1} nous dit que le morphisme
$$\chi_{mot} : K(\mathcal{V}(k)) \longrightarrow 
K(CMot(k))[\{(L^{i}-1)^{1}\}_{i>0}],$$
se factorise en un morphisme d'anneaux
$$\chi_{mot} : K(\mathcal{CH}^{sp}(k)) \longrightarrow 
K(CMot(k))[\{(L^{i}-1)^{1}\}_{i>0}].$$

\begin{df}\label{d11}
Pour un champ d'Artin sp\'ecial $F$, \emph{la caract\'eristique d'Euler motivique $F$} est
$\chi_{mot}([F])\in K(CMot(k))[\{(L^{i}-1)^{1}\}_{i>0}]$. Elle sera not\'ee
$\chi_{mot}(F)$.
\end{df}

En appliquant plusieurs foncteurs de r\'ealisations on trouve ainsi de nombreux invariants
pour les champs d'Artin sp\'eciaux. Par exemple, les nombres de Hodge de la d\'efinition
\ref{d6} peuvent \^etre extraits de la r\'ealisation de $\chi_{mot}(F)$ dans 
le groupe de Grothendieck des structures de Hodge pures. \\

\textit{Caract\'eristique d'Euler $l$-adique et formule des traces:}
Soit $k=\mathbb{F}_{q}$ un corps fini. Pour une vari\'et\'e $X$ sur $k$, on peut 
consid\'erer ses groupes de cohomologie \`a support compact
$H^{i}_{c}(\overline{X},\mathbb{Q}_{l})$, qui sont des $\mathbb{Q}_{l}$-espaces
vectoriels munis d'une action continue de $Gal(\overline{k}/k)\simeq \hat{\mathbb{Z}}$, c'est \`a dire
muni d'un automorphisme (on choisira le Frobenius g\'eom\'etrique). 
Les valeurs propres du Frobenius op\'erant sur $H^{i}_{c}(\overline{X},\mathbb{Q}_{l})$
sont des nombres de Weil relatif \`a $q$ (i.e. sont des nombres alg\'ebriques $\alpha$ tel que pour
tout plongement complexe on ait $|\alpha|=q^{\frac{n}{2}}$ pour un entier $n$). On notera $W(k)$ la cat\'egorie
ab\'elienne des $\mathbb{Q}_{l}$-espaces vectoriels de dimension finie
munis d'un automorphisme dont les valeurs propres sont des nombres de Weil (relatif \`a $q$). Cette cat\'egorie est 
une cat\'egorie tensorielle, et m\^eme Tannakienne. 
On consid\`ere alors $\Sigma_{i}(-1)^{i}[H^{i}_{c}(\overline{X},\mathbb{Q}_{l})]$ qui est un objet dans
l'anneau de Grothendieck $K(W(k))$. L'application
$$\mathcal{V}(k) \longrightarrow K(W(k))$$
qui \`a $X$ associe $\Sigma_{i}(-1)^{i}[H^{i}_{c}(\overline{X},\mathbb{Q}_{l})]$ 
se factorise en un morphisme d'anneaux
$$\chi_{l} : K(\mathcal{V}(k)) \longrightarrow K(W(k)).$$
De plus, $\chi_{l}(\mathbb{L})=\mathbb{Q}_{l}(1)$, o\`u $\mathbb{Q}_{l}(1)$ est de rang $1$ et l'action 
du Frobenius est donn\'ee par mutiplication par $q$ (il s'agit donc d'un \'el\'ement inversible).
Par le th\'eor\`eme \ref{t1} on trouve ainsi un morphisme d'anneaux
$$\chi_{l} : K(\mathcal{CH}^{sp}(k)) \longrightarrow K(W(k))[\{[\mathbb{Q}_{l}(i)-1]^{-1}\}_{i>0}],$$
o\`u $\mathbb{Q}_{l}(i)$ est de rang $1$ et l'action 
du Frobenius est donn\'ee par mutiplication par $q^{i}$. 

\begin{df}\label{d12}
Pour un champ d'Artin sp\'ecial $F$, \emph{la caract\'eristique d'Euler $l$-adique (\`a support
compact) de $F$} est
$\chi_{l}(F)\in K(W(k))[\{(\mathbb{Q}_{l}(i)-1)^{-1}\}_{i>0}].$
Elle sera not\'ee
$\chi_{l}(F)$.
\end{df}

A tout objet $V$ de $W(k)$, on peut associ\'e la trace du Frobenius g\'eom\'etrique, ce qui donne
un morphisme d'anneaux
$$Tr_{Fr} : K(W(k)) \longrightarrow \mathbb{Q}_{l},$$
qui est tel que $Tr_{Fr}([\mathbb{Q}_{l}(1)])=q$ et 
$Tr_{Fr}([\mathbb{Q}_{l}(i)-1])=q^{i}-1$ soient inversibles. 
On obtient ainsi un morphisme d'anneaux
$$Tr_{Fr}\circ \chi_{l} : K(\mathcal{CH}^{sp}(k)) \longrightarrow \mathbb{Q}_{l}.$$

La proposition suivante est une formule des traces pour les champs
d'Artin sp\'eciaux. C'est un cas particulier d'une formule beaucoup 
plus g\'en\'erale valable pour tout champ d'Artin fortement de 
pr\'esentation finie (voir \cite{tvv}). 

\begin{prop}\label{p5}
Pour tout $F \in \mathcal{CH}^{sp}(k)$, on a 
$$Tr_{Fr}(\chi_{l}(F)):=\mu(F),$$
o\`u $\mu(F)$ est le nombre de points rationels 
de $F$ d\'efini dans la proposition \ref{p4}.
\end{prop}

\textit{Preuve:} D'apr\`es le corollaire \ref{c3} il suffit de voir que
l'\'egalit\'e 
$$Tr_{Fr}(\chi_{l}(F)):=\mu(F)$$
est vraie lorsque $F$ est une vari\'et\'e. Mais ceci n'est autre que la
formule des traces de Grothendieck. 
\hfill $\Box$ \\

\end{document}